\documentstyle[12pt,leqno,a4wide,amssymb]{article}

\begin{document}

\title{\bf On cocycles with values in the group $SU(2)$
}
\author{Krzysztof Fr\c{a}czek}

\maketitle

\renewcommand{\thefootnote}{}

\newcommand{\tor}{{\Bbb T}}
\newcommand{\ad}{\,\mbox{Ad}}
\newcommand{\Pp}{{\cal P}}
\newcommand{\Hh}{{\cal H}}
\newcommand{\hsi}{\hat{\sigma}_f(n)}
\newcommand{\ep}{\varepsilon}
\newcommand{\pf}{{\bf Proof. }}
\newcommand{\su}{{\mathfrak s}{\mathfrak u}(2) }
\newcommand{\td}{\widetilde{\delta}}
\newcommand{\tp}{\widetilde{p} }
\newcommand{\tv}{\widetilde{\varphi}}
\newcommand{\ts}{\widetilde{\phi}}
\newcommand{\hp}{\widehat{\varphi}}
\newcommand{\jed}{\mbox{\boldmath$1$}}
\newcommand{\la}{\langle}
\newcommand{\pa}{\rangle}
\newtheorem{theo}{\indent Theorem}
[section]
\newtheorem{prop}[theo]{\indent Proposition}
\newtheorem{lem}[theo]{\indent Lemma}
\newtheorem{cor}[theo]{\indent Corollary}
\newtheorem{df}{\indent Definition}

\begin{abstract}
In this paper we introduce the notion of degree for $C^1$--cocycles over
irrational rotations on the circle with values in the group $SU(2)$. It is
shown that if a $C^1$--cocycle $\varphi:\tor\rightarrow SU(2)$ over an
irrational rotation by $\alpha$ has nonzero degree, then the skew product
\[\tor\times SU(2)\ni(x,g)\mapsto (x+\alpha,g\,\varphi(x))\in\tor\times SU(2)\] 
is not ergodic and the group of essential values of $\varphi$ is equal
to the maximal Abelian subgroup of $SU(2)$. Moreover, if $\varphi$ is of
class $C^2$ (with some additional assumptions) the Lebesgue component
in the spectrum of the skew product has countable multiplicity.  
Possible values of degree are discussed, too.
\end{abstract}

\section{Introduction}\indent

Assume that $T:(X,{\cal B},\lambda)\rightarrow(X,{\cal B},\lambda)$ is an
ergodic measure--preserving automorphism of standard Borel space.
Let $G$ be a compact Lie group,
$\mu$ its Haar measure. For a given measurable function
$\varphi:X\rightarrow G$ we study spectral properties of the 
measure--preserving automorphism of $X\times G$ (called skew product) 
defined by 
\[T_{\varphi}:(X\times G,\lambda\otimes\mu)\rightarrow(X\times
G,\lambda\otimes\mu),\;\; T_{\varphi}(x,g)=(Tx,g\,\varphi (x)).\]
A measurable function $\varphi:X\rightarrow G$ determines a measurable
cocycle over the automorphism $T$ given by
\[
\varphi^{(n)}(x)=\left\{
\begin{array}{ccl}
\varphi(x)\varphi(Tx)\ldots\varphi(T^{n-1}x) & \mbox{ for } & n>0\\
e & \mbox{ for } & n=0\\
(\varphi(T^{-1}x)\varphi(T^{-2}x)\ldots\varphi(T^{-n}x))^{-1} & \mbox{ for }
& n<0,
\end{array}
\right.
\] 
which we will identify with the function $\varphi$. Then
$T_{\varphi}^n(x,g)=(Tx,g\,\varphi^{(n)}(x))$ for any integer $n$.
Two cocycles $\varphi,\psi:X\rightarrow G$ are {\em cohomologous} if there
exists a measurable map $p:X\rightarrow G$ such that
\[\varphi(x)=p(x)^{-1}\,\psi(x)\,p(Tx).\]
In this case, $p$ will be called a transfer function. If $\varphi$ and
$\psi$ are cohomologous, then the map $(x,g)\mapsto(x,p(x)\,g)$
establishes a metrical isomorphism of $T_{\varphi}$ and $T_{\psi}$.

By $\tor$ we will mean the circle group $\{z\in{\Bbb C};|z|=1\}$
which most often will be treated as the group ${\Bbb R}/{\Bbb Z}$;
$\lambda$ will denote Lebesgue measure on $\tor$.
We will identify functions on $\tor$ with periodic of period 1 
functions  on ${\Bbb R}$. Assume that $\alpha\in\tor$ is irrational.
We will treat the case where $T$ is the ergodic rotation on $\tor$ given 
by $Tx=x+\alpha$.
 
In the case where $G$ is the circle and $\varphi$ is a smooth cocycle, 
spectral properties of $T_{\varphi}$ depend on the topological degree
$d(\varphi)$ of $\varphi$. For example,
in \cite{Iw-Le-Ru} A.\ Iwanik, M.\ Lema\'nczyk, D.\ Rudolph have proved
that if $\varphi$ is a $C^2$--cocycle with $d(\varphi)\neq 0$, then
$T_{\varphi}$ is ergodic and it has countable Lebesgue spectrum on the 
orthocomplement of the space of functions depending only on the
first variable.  On the other hand, in \cite{Ga-Le-Li} P.\ Gabriel, M.\
Lema\'nczyk, P.\ Liardet have proved that if $\varphi$ is absolutely
continuous with  $d(\varphi)=0$, then $T_{\varphi}$ has singular spectrum.

The aim of this paper is to find a spectral equivalent of topological
degree in case $G=SU(2)$.

\section{Degree of cocycle}\indent

In this section we introduce the notion of degree in case
$G=SU(2)$.  
For a given matrix
$A=[a_{ij}]_{i,j=1,\ldots,d}\in\mbox{M}_d({\Bbb C})$ define
$\|A\|=\sqrt{\frac{1}{d}\sum_{i,j=1}^d|a_{ij}|^2}$. Observe that if 
$A$ is an element of the Lie algebra $\su$, i.e.\ 
\[A=\left[
\begin{array}{cc}
ia & b+ic\\
-b+ic & -ia
\end{array}
\right],\]
where $a,b,c\in{\Bbb R}$, then $\|A\|=\sqrt{\det A}$.
Moreover, if $B$ is an element of the group $SU(2)$, i.e.\
\[A=\left[
\begin{array}{rr}
z_1 & z_2\\
-\overline{z_2} & \overline{z_1}
\end{array}
\right],\]
where $z_1,z_2\in{\Bbb C}$, $|z_1|^2+|z_2|^2=1$, then
$\ad_{B}A=BAB^{-1}\in\su$ and $\|\ad_BA\|=\|A\|$.  

Consider the scalar product of $\su$ given by
\[\la A,B\pa=-\frac{1}{8}\mbox{tr}(\mbox{ad}A\circ \mbox{ad}B).\]
Then $\|A\|=\sqrt{\la A,A\pa}$.
By $L^2(X,\su)$ we mean the space of all functions
$f:X\rightarrow\su$ such that
\[\|f\|_{L^2}=\sqrt{\int_{X}\|f(x)\|^2dx}<\infty.\]
For two $f_1,f_2\in L^2(X,\su)$ set
\[\la f_1,f_2\pa_{L^2}=\int_{X}\la f_1(x),f_2(x)\pa dx.\]
The space $L^2(X,\su)$ endowed with the above scalar product is a 
Hilbert space.

By $L^1(X,\su)$ we mean the space of all functions
$f:X\rightarrow\su$ such that
\[\|f\|_{L^1}=\int_{X}\|f(x)\|dx<\infty.\]
The space $L^1(X,\su)$ endowed with the norm $\|\:\|_{L^1}$ is a 
Banach space.

For a given measurable
cocycle $\varphi:\tor\rightarrow SU(2)$ consider the unitary operator
\begin{equation}\label{operat}
U:L^2(\tor,\su)\rightarrow L^2(\tor,\su),\;\;Uf(x)=\ad_{\varphi(x)}f(Tx).
\end{equation}
Then $U^nf(x)=\ad_{\varphi^{(n)}(x)}f(T^nx)$ for any
integer $n$.

\begin{lem}\label{pier}
There exists an operator $P:L^2(\tor,\su)\rightarrow L^2(\tor,\su)$ such
that
\[\frac{1}{n}\sum_{j=0}^{n-1} U^jf\rightarrow Pf\;\;\mbox{ in }\;\;
L^2(\tor,\su )\]
for any $f\in L^2(\tor,\su)$ and $U\circ P=P$. Moreover, $\|Pf\|$ is
constant $\lambda$-a.e.. 
\end{lem}

\pf 
First claim of the lemma follows from the von Neuman ergodic theorem. Since
$U\circ P=P$, we have $\ad_{\varphi(x)}Pf(Tx)=Pf(x)$,
for $\lambda$-a.e.\ $x\in\tor$. It follows
that $\|Pf(Tx)\|=\|Pf(x)\|$, for $\lambda$-a.e.\ $x\in\tor$. Hence
$\|Pf(x)\|=c$, for $\lambda$-a.e.\ $x\in\tor$, by the ergodicity of $T$.
$\blacksquare$ 

\begin{lem}\label{drugi}
For every $f\in L^2(\tor,\su)$ the sequence
$\frac{1}{n}\sum_{j=0}^{n-1} U^jf$
converges $\lambda$--almost everywhere. 
\end{lem}

\pf Let $\tilde{f}\in L^2(\tor\times SU(2),\su)$ be given by
$\tilde{f}(x,g)=\ad_{g}f(x)$. Then
\[\tilde{f}(T_{\varphi}^n(x,g))=\ad_{g}(U^nf(x))\]
for any integer $n$.
By the Birkhoff ergodic theorem, 
the sequence
\[\frac{1}{n}\sum_{j=0}^{n-1}\tilde{f}(T_{\varphi}^n(x,g))=
\ad_{g}(\frac{1}{n}\sum_{j=0}^{n-1} U^jf(x))\]
converges for $\lambda\otimes\mu$--a.e.\ $(x,g)\in\tor\times SU(2)$. Hence
there exists $g\in SU(2)$ such that $\ad_g(\frac{1}{n}\sum_{j=0}^{n-1}
U^jf(x))$ converges for $\lambda$-a.e.\ $x\in\tor$, and the proof is
complete. $\blacksquare$ 
\vspace{2ex}

Recall that, if a function $\varphi:\tor\rightarrow SU(2)$ is of class
$C^1$, then $D\varphi(x)\varphi(x)^{-1}\in\su$ for every $x\in\tor$.

\begin{lem}\label{defin}
For every $C^1$--cocycle $\varphi:\tor\rightarrow SU(2)$, there
exists  $\psi\in L^2(\tor,\su)$ such that
\[\frac{1}{n}D\varphi^{(n)}(\varphi^{(n)})^{-1}\rightarrow\psi\;\;
\mbox{ in $L^2(\tor,\su)$ and $\lambda$--almost everywhere.}\]
Moreover, $\|\psi\|$ is a constant function $\lambda$--a.e.\ and
$\varphi(x)\psi(Tx)\varphi(x)^{-1}=\psi(x)$ for $\lambda$-a.e.\ $x\in\tor$.
\end{lem}

\pf
Since
\[D\varphi^{(n)}(x)=\sum_{j=0}^{n-1}\varphi(x)\ldots\varphi(T^{j-1}x)D\varphi(T^{j}x)\varphi(T^{j+1}x)\ldots\varphi(T^{n-1}x),\]
we have
\begin{eqnarray*}
D\varphi^{(n)}(x)(\varphi^{(n)}(x))^{-1} & = & 
\sum_{j=0}^{n-1}\varphi(x)\ldots\varphi(T^{j-1}x)D\varphi(T^{j}x)\varphi(T^{j}x)^{-1}\varphi(T^{j-1}x)^{-1}\ldots\varphi(x)^{-1}\\
 & = & 
\sum_{j=0}^{n-1}\varphi^{(j)}(x)D\varphi(T^{j}x)\varphi(T^{j}x)^{-1}(\varphi^{(j)}(x))^{-1}\\
 & =&
\sum_{j=0}^{n-1}U^j(D\varphi\:\varphi^{-1})(x),
\end{eqnarray*}
where $U$ is the unitary operator
given by (\ref{operat}). Applying Lemmas~\ref{pier} and \ref{drugi}, we
get  $\psi=P(D\varphi\:\varphi^{-1})$, and the 
proof is complete. $\blacksquare$

\begin{df}
{\em The number $\|\psi\|$ will be called the} degree {\em  of the cocycle
$\varphi$ and denoted by $d(\varphi)$.}
\end{df}

Lemma~\ref{defin} shows that 
\[\frac{1}{n}\|D\varphi^{(n)}(\varphi^{(n)})^{-1}\|_{L^1}\rightarrow
d(\varphi).\] 
On the other hand, $\|D\varphi^{(n)}(\varphi^{(n)})^{-1}\|_{L^1}$ is the
length of  the
curve $\varphi^{(n)}$. Geometrically speaking, the degree of $\varphi$
is the limit of length($\varphi^{(n)})/n$.

A measurable cocycle $\delta:\tor\rightarrow SU(2)$ is called {\em diagonal} if 
there exists a measurable function $\gamma:\tor\rightarrow\tor$ such that
\[\delta(x)=\left[
\begin{array}{cc}
\gamma(x) & 0\\
0 & \overline{\gamma(x)}
\end{array}\right].\]

\begin{theo}\label{nieerg}
Suppose that $\varphi:\tor\rightarrow SU(2)$ is a $C^1$--cocycle  with
$d(\varphi)\neq 0$. Then $\varphi$ is cohomologous to a diagonal
cocycle.
\end{theo}

\pf
For every nonzero $A\in\su$ there exists $B_A\in SU(2)$ such that
\[B_AA(B_A)^{-1}=\left[
\begin{array}{cc}
i\|A\| & 0\\
0 & -i\|A\|
\end{array}\right].\]
Indeed, if
$A=\left[
\begin{array}{cc}
ia & b+ic\\
-b+ic & -ia
\end{array}
\right],$
then we can take
\[B_A=\left\{
\begin{array}{ccl}
\left[
\begin{array}{cc}
-i\sqrt{\frac{\|A\|+a}{2\|A\|}}\frac{b+ic}{|b+ic|} & 
-\sqrt{\frac{\|A\|-a}{2\|A\|}}\\
\sqrt{\frac{\|A\|-a}{2\|A\|}} & 
i\sqrt{\frac{\|A\|+a}{2\|A\|}}\frac{b-ic}{|b+ic|}
\end{array}\right]
& \mbox{if} & |a|\neq\|A\|\\
\left[
\begin{array}{cr}
0 & -1\\
1 & 0
\end{array}\right] 
& \mbox{if} & a=-\|A\|\\
\left[
\begin{array}{cr}
1 & 0\\
0 & 1
\end{array}\right] 
& \mbox{if} & a=\|A\|.
\end{array}\right.\]
Set $p(x)=B_{\psi(x)}$. Then $p:\tor\rightarrow SU(2)$ is a measurable
function and
\[\psi(x)=p(x)^{-1}\left[
\begin{array}{cc}
i\:d(\varphi) & 0\\
0 & -i\:d(\varphi)
\end{array}
\right]p(x).\]
Since $\varphi(x)\psi(Tx)\varphi(x)^{-1}=\psi(x)$, we have
\[\varphi(x)p(Tx)^{-1}\left[
\begin{array}{cc}
i\:d(\varphi) & 0\\
0 & -i\:d(\varphi)
\end{array}
\right]p(Tx)\varphi(x)^{-1}
=p(x)^{-1}\left[
\begin{array}{cc}
i\:d(\varphi) & 0\\
0 & -i\:d(\varphi)
\end{array}
\right]p(x).\]
Hence
\[p(x)\varphi(x)p(Tx)^{-1}\left[
\begin{array}{cc}
i\:d(\varphi) & 0\\
0 & -i\:d(\varphi)
\end{array}
\right]
=\left[
\begin{array}{cc}
i\:d(\varphi) & 0\\
0 & -i\:d(\varphi)
\end{array}
\right]p(x)\varphi(x)p(Tx)^{-1}.\]
Since $d(\varphi)\neq 0$, we see that the cocycle
$\delta:\tor\rightarrow SU(2)$ defined by
$\delta(x)=p(x)\varphi(x)p(Tx)^{-1}$ is diagonal. $\blacksquare$

\vspace{2ex}

For a given $C^1$--cocycle $\varphi:\tor\rightarrow SU(2)$ with nonzero
degree let $\gamma=\gamma(\varphi):\tor\rightarrow\tor$ be a measurable
cocycle such that the cocycles $\varphi$ and $\left[
\begin{array}{cc}
\gamma & 0\\
0 & (\gamma)^{-1}
\end{array}\right]$
are cohomologous. It is easy to check that the
choice of $\gamma$ is unique up to a measurable cohomology with values in
the circle and inverse.  

Theorem~\ref{nieerg} shows that if $d(\varphi)\neq 0$, then the skew product
$T_{\varphi}$ is metrically isomorphic to a skew product of an irrational 
rotation on the circle and a diagonal cocycle. It follows that
$T_{\varphi}$ is not ergodic. However, in the next sections we show that
if $d(\varphi)\neq 0$, then $\varphi$ is not cohomologous to a constant
cocycle. Moreover, the skew product
$T_{\gamma}:\tor\times\tor\rightarrow\tor\times\tor$ is ergodic and it is
mixing on the  orthocomplement of the space of functions
depending only on the first variable. We prove also that 
(with some additional assumptions on $\varphi$) the
Lebesgue component in the spectrum of $T_{\gamma}$ has countable multiplicity.
It follows that if $d(\varphi)\neq 0$, then:
\begin{itemize}
\item all ergodic components of $T_{\varphi}$ are metrically isomorphic to
$T_{\gamma}$, 
\item the spectrum of $T_{\varphi}$ consists of two parts: discrete and
 mixing,
\item (with some additional assumptions on $\varphi$) the
Lebesgue component in the spectrum of $T_{\varphi}$ has countable
multiplicity.  
\end{itemize}

In case $G=\tor$ the topological degree of each $C^1$--cocycle is an
integer number. An important question is: what can one say on values of
degree in case $G=SU(2)$?  
If a cocycle $\varphi$ is cohomologous to a diagonal cocycle via a smooth
transfer function, then $d(\varphi)\in 2\pi {\Bbb N}_0= 
2\pi ({\Bbb N}\cup \{0\})$.

We call a function $f:\tor\rightarrow SU(2)$ absolutely continuous
if $f_{ij}:\tor\rightarrow{\Bbb C}$ is absolutely continuous for $i,j=1,2$.
Suppose that $\varphi$ is cohomologous to a diagonal cocycle via an
absolutely continuous transfer function. Then $\varphi$ can be represented as 
$\varphi(x)=p(x)^{-1}\delta(x)p(Tx)$, where 
$\delta,p:\tor\rightarrow SU(2)$ are absolutely continuous and 
$\delta$ is diagonal. Since
$\varphi^{(n)}(x)=p(x)^{-1}\delta^{(n)}(x)p(T^nx)$, we have
\begin{eqnarray*}
\frac{1}{n}D\varphi^{(n)}(x)(\varphi^{(n)}(x))^{-1} & = &
\frac{1}{n}(-p(x)^{-1}Dp(x)+
\varphi^{(n)}(x)p(T^nx)^{-1}Dp(T^nx)(\varphi^{(n)}(x))^{-1}\\ & &+
p(x)^{-1}D\delta^{(n)}(x)(\delta^{(n)}(x))^{-1}p(x)).
\end{eqnarray*} 
On the other hand, $\delta(x)=\left[
\begin{array}{cc}
\gamma(x) & 0\\
0 & \overline{\gamma(x)}
\end{array}
\right],$
where $\gamma:\tor\rightarrow\tor$ is an absolutely continuous cocycle of
the
form $\gamma(x)=\exp 2\pi i(\widetilde{\gamma}(x)+kx)$, where $k$ is the
topological 
degree of $\gamma$ and $\widetilde{\gamma}:\tor\rightarrow{\Bbb R}$ is
an absolutely continuous function. 
Then
\[\frac{1}{n}D\gamma^{(n)}(x)(\gamma^{(n)}(x))^{-1}=
2\pi i(\frac{1}{n}\sum_{j=0}^{n-1}D\widetilde{\gamma}(T^jx)+k)
\rightarrow 2\pi ik\] 
in $L^1(\tor,{\Bbb R})$, by the Birkhoff ergodic theorem. It follows that
\[\frac{1}{n}D\varphi^{(n)}(x)(\varphi^{(n)}(x))^{-1}\rightarrow
p(x)^{-1}
\left[
\begin{array}{cc}
2\pi ik & 0\\
0 & -2\pi ik
\end{array}
\right]
p(x)\]
in $L^1(\tor,\su)$.
Hence $d(\varphi)=2\pi|d(\gamma)|\in 2\pi {\Bbb N}_0$.

In Section~\ref{integ} it is shown that if $\alpha$ is the golden ration,
then the degree of every $C^2$--cocycle belongs to $2\pi {\Bbb N}_0$, too. 

\section{Notation and facts from spectral theory}\indent

Let $U$ be a unitary operator on a separable Hilbert space
${\cal H}$.
By the {\em cyclic space} generated by $f\in {\cal H}$ we mean  the space 
${\Bbb Z}(f)=\mbox{span}\{U^nf;n\in{\Bbb Z}\}$.
By the {\em spectral measure} $\sigma_f$ of $f$ we mean a Borel
measure on $\tor$ determined by the equalities
\[\hat{\sigma}_f(n)=\int_{\tor} e^{2\pi inx}d\sigma_f(x)=(U^nf,f)\]
for $n\in{\Bbb Z}$. Recall that
there exists a sequence $\{f_n\}_{n\in{\Bbb N}}$ in ${\cal H}$ such that
\begin{equation}\label{spectral}
\begin{array}{ccc}
{\cal H}=\bigoplus_{n=1}^{\infty} {\Bbb Z}(f_n) & and &
\sigma_{f_1}\gg\sigma_{f_2}\ldots 
\end{array}.
\end{equation}
Moreover, for any sequence $\{f_n'\}_{n\in{\Bbb N}}$ in ${\cal H}$ satisfying
(\ref{spectral}) we have
$\sigma_{f_1}\equiv\sigma_{f_1'},\sigma_{f_2}\equiv\sigma_{f_2'},\ldots$ .
The above decompositions of $\Hh$ are called spectral decompositions
of $U$.

The spectral type of $\sigma_{f_1}$ (the equivalence class of
measures) will be called the {\em maximal spectral type} of
$U$. We say that $U$ has {\em Lebesgue (continuous singular,
discrete) spectrum} if $\sigma_{f_1}$ is equivalent to Lebesgue
(continuous singular, discrete) measure on the circle.
An operator $U$ is called  mixing if 
\[\hat{\sigma}_f(n)=(U^nf,f)\rightarrow 0\]
for any $f\in{\cal H}$.
We say that the Lebesgue component in the spectrum of $U$ has countable
multiplicity if $\lambda\ll\sigma_{f_n}$ for every natural $n$ or
equivalently if there exists a sequence $\{g_n\}_{n\in{\Bbb N}}$ in ${\cal
H}$ such that the cyclic spaces ${\Bbb Z}(g_n)$ are pairwise
orthogonal and $\sigma_{g_n}\equiv\lambda$ for every natural $n$.

For a skew product $T_{\varphi}$ consider its Koopman operator
\[U_{T_{\varphi}}:L^2(\tor\times G,\lambda\otimes\mu)\rightarrow
L^2(\tor\times G,\lambda\otimes\mu),\;\; 
U{T_{\varphi}}f(x,g)=f(Tx,g\,\varphi (x)).\]
Denote by $\widehat{G}$ the set of all equivalence classes of unitary
irreducible representations of the group $G$. For any unitary irreducible 
representation $\Pi:G\rightarrow{\cal U}({\cal H}_{\Pi})$
by $\{\Pi_{ij}\}_{i,j=1}^{d_{\Pi}}$ we mean the matrix elements of $\Pi$,
where $d_{\Pi}=\dim {\cal H}_{\Pi}$.
Let us decompose
\[L^2(\tor\times G)=\bigoplus_{\Pi\in\widehat{G}}\bigoplus_{i=1}^{d_{\pi}}
\Hh_i^{\Pi},\] 
where
\[\Hh_i^{\Pi}=\{\sum_{j=1}^{d_{\pi}}\Pi_{ij}(g)f_j(x);f_j\in
L^2(\tor,\lambda),j=1,\ldots,d_{\Pi}\} \simeq
\overbrace{L^2(\tor,\lambda)\oplus\ldots\oplus L^2(\tor,\lambda)}^{d_{\Pi}}.\]
Observe that $\Hh_i^{\Pi}$ is a closed $U_{T_{\varphi}}$--invariant
subspace of $L^2(\tor\times G)$ and
\[U^n_{T_{\varphi}}(\sum_{j=1}^{d_{\pi}}\Pi_{ij}(g)f_j(x))=
\sum_{j,k=1}^{d_{\pi}}\Pi_{ik}(g)\Pi_{kj}(\varphi^{(n)}(x))f_j(T^nx).\]
Consider the unitary operator
$M_i^{\Pi}:\Hh_i^{\Pi}\rightarrow\Hh_i^{\Pi}$ given by 
\[M_i^{\Pi}(\sum_{j=1}^{d_{\pi}}\Pi_{ij}(g)f_j(x))=
\sum_{j=1}^{d_{\pi}}e^{2\pi ix}\Pi_{ij}(g)f_j(x).\]
Then
\begin{equation}\label{komut}
U_{T_{\varphi}}^nM_i^{\Pi}f=e^{2\pi in\alpha}M_i^{\Pi}U_{T_{\varphi}}^nf
\end{equation}
for any $f\in\Hh_i^{\Pi}$.
It follows that 
\[\int_{\tor} e^{2\pi inx}d\sigma_{M_i^{\Pi}f}(x)=
(U^n_{T_{\varphi}}M_i^{\Pi}f,M_i^{\Pi}f)=
e^{2\pi in\alpha}(U^n_{T_{\varphi}}f,f)=
\int_{\tor} e^{2\pi inx}d(T^*\sigma_{f})(x)\]
for any $f\in\Hh_i^{\Pi}$.
Hence $\sigma_{M_i^{\Pi}f}=T^*\sigma_{f}$.

\begin{lem}\label{unif}
For every $\Pi\in\widehat{G}$ and $i=1,\ldots,d_{\pi}$ if the operator 
$U_{T_{\varphi}}:\Hh_i^{\Pi}\rightarrow\Hh_i^{\Pi}$ has
absolutely continuous spectrum, then it has Lebesgue spectrum of uniform
multiplicity.
\end{lem}

\pf Let $\Hh_i^{\Pi}=\bigoplus_{n=1}^{\infty} {\Bbb Z}(f_n)$ be a spectral
decomposition. Then
\[\Hh_i^{\Pi}=(M_i^{\Pi})^m\Hh_i^{\Pi}=
\bigoplus_{n=1}^{\infty} {\Bbb Z}((M_i^{\Pi})^mf_n)\]
is a spectral decomposition for any integer $m$. Therefore
$\sigma_{f_n}\equiv\sigma_{(M_i^{\Pi})^mf_n}\ll\lambda$ for every natural
$n$ and integer $m$.
Suppose that there exists a Borel set $A\subset\tor$ such that
$\sigma_{f_n}(A)=0$ and $\lambda(A)>0$. Then
\[\sigma_{f_n}(\bigcup_{m\in{\Bbb Z}}T^mA)=0\;\;\mbox{ and }\;\;
\lambda(\bigcup_{m\in{\Bbb Z}}T^mA)=1,\]
by the ergodicity of $T$.
It follows that $\sigma_{f_n}\equiv\lambda$ or $\sigma_{f_n}=0$ for every
natural $n$. $\blacksquare$

\begin{lem}
If
\[\sum_{n\in{\Bbb Z}}|\int_{\tor}\Pi_{jj}(\varphi^{(n)}(x))dx|^2<\infty\]
for $j=1,\ldots,d_{\Pi}$, then $U_{T_{\varphi}}$ has Lebesgue spectrum of
uniform multiplicity on $\Hh_{i}^{\Pi}$ for $i=1,\ldots,d_{\Pi}$.
\end{lem}

\pf Fix $1\leq i\leq d_{\Pi}$. Note that
\[\la U^n_{T_{\varphi}}\Pi_{ij},\Pi_{ij}\pa=
\sum_{k=1}^{d_{\Pi}}\int_{\tor}\int_{G}
\la\Pi_{ik}(g)\Pi_{kj}(\varphi^{(n)}(x)),\Pi_{ij}(g)\pa dgdx=
\frac{1}{d_{\Pi}}\int_{\tor}\Pi_{jj}(\varphi^{(n)}(x))dx.\]
Since
\[\sum_{n\in{\Bbb Z}}|\la U^n_{T_{\varphi}}\Pi_{ij},\Pi_{ij}\pa|^2<\infty,\]
we have $\sigma_{\Pi_{ij}}\ll\lambda$ for $j=1,\ldots,d_{\Pi}$.
From (\ref{komut}) we  get $\sigma_{(M_i^{\Pi})^m\Pi_{ij}}\ll\lambda$ for
any integer $m$.  Since $\{f\in\Hh_i^{\Pi};\sigma_{f}\ll\lambda\}$ is a
closed linear subspace of $L^2(\tor\times G)$ and the set 
$\{(M_i^{\Pi})^m\Pi_{ij};j=1,\ldots,d_{\Pi},m\in{\Bbb Z}\}$  generates the
space $\Hh_i^{\Pi}$, it follows that $U_{T_{\varphi}}$ has absolutely
continuous spectrum on  $\Hh_i^{\Pi}$. By Lemma~\ref{unif}, $U_{T_{\varphi}}$
has Lebesgue spectrum of uniform multiplicity on $\Hh_i^{\Pi}$. $\blacksquare$

\begin{cor}\label{kica}
For every $\Pi\in\widehat{G}$, if
\[\sum_{n\in{\Bbb Z}}\|\int_{\tor}\Pi(\varphi^{(n)}(x))dx\|^2<\infty,\]
then $U_{T_{\varphi}}$ has Lebesgue spectrum of
uniform multiplicity on $\bigoplus_{i=1}^{d_{\Pi}}\Hh_{i}^{\Pi}$.
$\blacksquare$ 
\end{cor}
 
Similarly one can prove the following result.

\begin{theo}\label{mix}
For every $\Pi\in\widehat{G}$, if
\[\lim_{n\rightarrow\infty}\int_{\tor}\Pi(\varphi^{(n)}(x))dx=0,\]
then $U_{T_{\varphi}}$ is mixing on
$\bigoplus_{i=1}^{d_{\Pi}}\Hh_{i}^{\Pi}$. $\blacksquare$ 
\end{theo}

\section{Representations of $SU(2)$}\indent

In this section some basic information about the theory of
representations of the group $SU(2)$ are presented. By $\Pp_k$ we mean the
linear space of all homogeneous  
polynomials of degree $k\in{\Bbb N}_0$ in two variables $u$ and $v$.
Denote by $\Pi_k$ the representation of the group $SU(2)$ in
$\Pp_k$ given by 
\[\left[ \Pi_k\left(\left[
\begin{array}{rr}
z_1 & z_2\\
-\overline{z_2} & \overline{z_1}
\end{array}
\right]\right)f\right](u,v)=f(z_1u-\overline{z_2}v,z_2u+\overline{z_1}v).\]
Of course, all $\Pi_k$ are unitary (under an appropriate
inner product on $\Pp_k$) and the family $\{\Pi_0,\Pi_1,\Pi_2,\ldots\}$ is
a complete family of continuous unitary irreducible representations of
$SU(2)$. In the Lie algebra $\su$, we choose the following basis:
\[
\begin{array}{ccc}
h=\left[
\begin{array}{rr}
1 & 0\\
0 & -1
\end{array}
\right], &
e=\left[
\begin{array}{rr}
0 & 1\\
0 & 0
\end{array}
\right], &
f=\left[
\begin{array}{rr}
0 & 0\\
1 & 0
\end{array}
\right].
\end{array}
\]
Let $V_k$ be a $k+1$-dimension linear space. For every natural $k$ there
exists a basis $v_0,\ldots,v_k$ of $V_k$ such that the corresponding
representation $\Pi_k^*$ of $\su$ in $V_k$ has the following form:
\begin{eqnarray*}
\Pi_k^*(e)v_i & = & i(k-i+1)v_{i-1}\\
\Pi_k^*(f)v_i & = & v_{i+1}\\
\Pi_k^*(h)v_i & = & (k-2i)v_{i}
\end{eqnarray*}
for $i=0,\ldots,k$. Then 
\begin{equation}\label{taki}
\|A\|\leq\|\Pi_k^*(A)\|\leq k^2\|A\|
\end{equation}
for any $A\in\su$.
\begin{lem}\label{kocie}
\[\det\Pi_{2k-1}^*(A)=((2k-1)!!)^2(\det A)^k\]
for any $A\in\su$ and $k\in{\Bbb N}$.
\end{lem}

\pf
For every $A\in\su$ there exists $B\in SU(2)$ and $d\in{\Bbb R}$ such that 
$A=\ad_{B}\left[
\begin{array}{cc}
id & 0\\
0 & -id
\end{array}
\right]$. Then
\[\Pi_{2k-1}^*(A)=\Pi_{2k-1}^*(\ad_B\left[
\begin{array}{cc}
id & 0\\
0 & -id
\end{array}
\right])=
\ad_{\Pi_{2k-1}(B)}\Pi_{2k-1}^*(\left[
\begin{array}{cc}
id & 0\\
0 & -id
\end{array}
\right]).\]
It follows that
\[\det\Pi_{2k-1}^*(A)=\det\Pi_{2k-1}^*(\left[
\begin{array}{cc}
id & 0\\
0 & -id
\end{array}
\right])=((2k-1)!!)^2d^{2k}=((2k-1)!!)^2(\det A)^k.\;\;\blacksquare\]

\begin{lem}\label{kickic}
For any nonzero $A\in\su$ the matrix $\Pi_{2k-1}^*(A)$ is invertible.
Moreover, for every natural $k$ there exists a real constant $K_k>0$ such
that 
\[\|\Pi_{2k-1}^*(A)^{-1}\|\leq K_k\|A\|^{-1}\]
for every nonzero $A\in\su$.
\end{lem}

\pf First claim of the lemma follows from Lemma~\ref{kocie}.
Set $C=\Pi_{2k-1}^*(A)$. Then 
\[|[C]_{ij}|\leq (2k)^{4k}(2k-1)!\|A\|^{2k-1}\]
for $i,j=1,\ldots,2k$. It follows that
\[|(C^{-1})_{ij}|=\frac{|[C]_{ij}|}{\det\Pi_{2k-1}^*(A)}\leq
\frac{(2k)^{4k}(2k-1)!\|A\|^{2k-1}}{((2k-1)!!)^2\|A\|^{2k}}=
\frac{(2k)^{4k}(2k-1)!}{((2k-1)!!)^2}\|A\|^{-1}.\]
Hence
\[\|C^{-1}\|\leq \frac{(2k)^{4k+1}(2k-1)!}{((2k-1)!!)^2}\|A\|^{-1}.
\;\;\blacksquare\]

\section{Ergodicity and mixing of $T_{\gamma}$}
\begin{lem}
Suppose that $\{f_n\}_{n\in{\Bbb N}}$ is a sequence in $L^2(\tor,{\Bbb
C})$ such that $\int_0^xf_n(y)dy\rightarrow 0$ for any $x\in\tor$. Let
$g:\tor\rightarrow{\Bbb C}$ be a bounded measurable function.
Then
\[\lim_{n\rightarrow\infty}
\int_{\tor} f_n(y)g(T^ny)dy= 0\;\;\mbox{ and }\;\;
\lim_{n\rightarrow\infty}\int_0^xf_n(y)g(y)dy=0\]
for any $x\in\tor$.
\end{lem}

\pf
By assumption, the sequence $\{f_n\}_{n\in{\Bbb N}}$ tends to zero in the weak
topology in $L^2(\tor,{\Bbb C})$, which implies immediately the second
claim of the lemma. 
Since $\{f_n\}_{n\in{\Bbb N}}$ converges weakly to zero, for every integer
$m$ we have
\[\lim_{n\rightarrow\infty}\int_{\tor}f_n(T^{-n}y)\exp 2\pi imy\,dy=
\lim_{n\rightarrow\infty}\int_{\tor}f_n(y)\exp 2\pi im(y+n\alpha)\,dy=0.\]
It follows that the sequence $\{f_n\circ T^{-n}\}_{n\in{\Bbb N}}$
converges weakly to zero. Therefore
\[\int_{\tor} f_n(y)g(T^ny)dy=\int_{\tor} f_n(T^{-n}y)g(y)dy=0.
\;\blacksquare\]

This gives immediately the following conclusion.

\begin{cor}\label{cormix}
Suppose that $\{f_n\}_{n\in{\Bbb N}}$ is a sequence in 
$L^2(\tor,M_k({\Bbb C}))$ ($k$ is a natural number) such
that $\int_0^xf_n(y)dy\rightarrow 0$ for any $x\in\tor$. Let
$g:\tor\rightarrow M_k({\Bbb C})$ be a bounded measurable function.
Then
\[\lim_{n\rightarrow\infty}
\int_{\tor} f_n(y)g(T^ny)dy= 0\;\;\mbox{ and }\;\;
\lim_{n\rightarrow\infty}\int_0^xf_n(y)g(y)dy=0\]
for any $x\in\tor$. $\blacksquare$
\end{cor}

\begin{theo}
Let $\varphi:\tor\rightarrow SU(2)$ be a $C^1$--cocycle with nonzero degree.
Then the skew product
$T_{\gamma(\varphi)}:\tor\times\tor\rightarrow\tor\times\tor$ is ergodic 
and it is mixing on the  orthocomplement of the space of functions
depending only on the first variable.
\end{theo}

\pf By Theorem~\ref{mix}, it suffices to show that 
\[\lim_{n\rightarrow\infty}\int_{\tor}(\gamma^{(n)}(x))^kdx=0\]
for every nonzero integer $k$.
Fix $k\in{\Bbb N}$. 
Denote by $\psi:\tor\rightarrow\su$ the limit (in $L^2(\tor,\su)$) of the
sequence
$\{\frac{1}{n}D\varphi^{(n)}(\varphi^{(n)})^{-1}d\lambda\}_{n\in{\Bbb N}}$.
Let $p:\tor\rightarrow SU(2)$ be a measurable function
such that
\[\left[
\begin{array}{cc}
 \gamma(x) & 0 \\
0 & \overline{\gamma(x)}
\end{array}
\right]=
p(x)\varphi(x)p(Tx)^{-1}
\;\;\mbox{ and }\;\;
\ad_{p(x)}(\psi(x))=
\left[
\begin{array}{lr}
id & 0 \\
0 & -id
\end{array}
\right],\]
where $d$ is the degree of $\varphi$ (see the proof of
Theorem~\ref{nieerg}). Then
\begin{equation}\label{mat1}
\left[
\begin{array}{clcrc}
 (\gamma^{(n)})^k & & & & \\
 & (\gamma^{(n)})^{k-2} & & 0 & \\
 &  & \ddots &  & \\
 & 0 & & (\gamma^{(n)})^{-k+2} & \\
 & & & & (\gamma^{(n)})^{-k}
\end{array}
\right]
=
\Pi_k(p)\Pi_k(\varphi^{(n)})\Pi_k(p\circ T^n)^{-1}
\end{equation}
for any natural $n$ and
\begin{eqnarray}\label{mat2}
\ad_{\Pi_k(p(x))}\Pi_k^*(\psi(x)) & = & 
\Pi_k^*(\ad_{p(x)}\psi(x))\\
 & = & \nonumber
\Pi_k^*
(\left[
\begin{array}{cr}
 id & 0 \\
 0 & -id
\end{array}
\right])
\\
 & = & \nonumber
\left[
\begin{array}{clcrc}
 kid & & & & \\
 & (k-2)id & & 0 & \\
 &  & \ddots &  & \\
 & 0 & & (-k+2)id & \\
 & & & & -kid
\end{array}
\right].
\end{eqnarray}
Recall that for any differentiable function $\xi:\tor\rightarrow SU(2)$ and
for any representation $\Pi$ of $SU(2)$ we have
\[D(\Pi\xi(x))(\Pi\xi(x))^{-1}=\Pi^*(D\xi(x)\xi(x)^{-1}).\]
Therefore
\begin{eqnarray*}
\int_0^x\frac{1}{n}\Pi_k^*(D\varphi^{(n)}(y)(\varphi^{(n)}(y))^{-1})
\Pi_k(\varphi^{(n)}(y))dy & = &
\int_0^x\frac{1}{n}D(\Pi_k\varphi^{(n)}(y))dy\\
& = &
\frac{1}{n}(\Pi_k(\varphi^{(n)}(x))-\Pi_k(\varphi^{(n)}(0)))
\end{eqnarray*}
tends to zero for any $x\in\tor$. Since 
\[\frac{1}{n}\Pi_k^*(D\varphi^{(n)}(\varphi^{(n)})^{-1})\rightarrow 
\Pi_k^*\psi\]
in $L^2(\tor,M_{k+1}({\Bbb C}))$, it follows that
\[\int_0^x\Pi_k^*(\psi(y))
\Pi_k(\varphi^{(n)}(y))dy\rightarrow 0 \]
 for any $x\in\tor$.
By Corollary~\ref{cormix},
\[\int_{\tor}\Pi_k(p(y))\Pi_k^*(\psi(y))
\Pi_k(\varphi^{(n)}(y))\Pi_k(p(T^ny))^{-1}dy\rightarrow 0.\]
On the other hand,
\[\Pi_k(p(y))\Pi_k^*(\psi(y))
\Pi_k(\varphi^{(n)}(y))\Pi_k(p(T^ny))^{-1}=
\left[\begin{array}{ccc}
ikd(\gamma^{(n)}(y))^k & & 0\\
 & \ddots & \\
 0 & & -ikd(\gamma^{(n)}(y))^{-k}
 \end{array}\right],\]
by (\ref{mat1}) and (\ref{mat2}).
Therefore 
\[\lim_{n\rightarrow\infty}\int_{\tor}(\gamma^{(n)}(y))^mdy=0\]
for any nonzero $m\in\{-k,-k+2,\ldots,k-2,k\}$, which completes the proof.
$\blacksquare$ 

\section{Spectral analysis of cocycles with nonzero degree}\indent

In this section it is shown that for every cocycle 
$\varphi:\tor\rightarrow SU(2)$ if $d(\varphi)\neq 0$ and it satisfies some
additional assumptions, then the Lebesgue component in the spectrum of
$T_{\varphi}$ has countable multiplicity.

Now we introduce a notation, which is necessary to prove the main theory.  
Let $f,g:\tor\rightarrow M_{k}({\Bbb C})$ be functions of bounded variation
(i.e.\ $f_{ij},g_{ij}:\tor\rightarrow{\Bbb C}$ have bounded variation for
$i,j=1,\ldots,k$) 
and let one of them be continuous. We will use the symbol
$\int_{\tor}f\,dg$ to denote the $k\times k$-matrix given by
\[(\int_{\tor}f\,dg)_{ij}=\sum_{l=1}^k\int_{\tor}f_{il}\,dg_{lj}\]
for $i,j=1,\ldots,d$. It is clear that if $g$ is absolutely continuous, then
\begin{equation}\label{cal1}
\int_{\tor}f\:dg=\int_{\tor}f(x)\:Dg(x)dx.
\end{equation}
Moreover, applying integration by parts, we have
\begin{equation}\label{cal2}
\int_{\tor}f\:dg=-(\int_{\tor}g^T\:df^T)^T.
\end{equation}

\begin{theo}\label{glowne1}
Let $\varphi:\tor\rightarrow SU(2)$ be a $C^2$--cocycle  with
$d(\varphi)\neq 0$. Suppose that the sequence
$\{\frac{1}{n}D\varphi^{(n)}(\varphi^{(n)})^{-1}\}_{n\in{\Bbb N}}$
is uniformly convergent and
$\{D(\frac{1}{n}D\varphi^{(n)}(\varphi^{(n)})^{-1})\}_{n\in{\Bbb N}}$
is bounded in $L^1(\tor,\su)$. Then the Lebesgue component in the spectrum of
$T_{\varphi}$ has countable multiplicity.
Moreover, the Lebesgue component in the spectrum of
$T_{\gamma(\varphi)}$ has countable multiplicity, too.
\end{theo}

\pf
First observe that it suffices to show that for every natural $k$
there exists a real constant $C_k>0$ such that
\begin{equation}\label{lebe}
\|\int_{\tor}\Pi_{2k-1}(\varphi^{(n)}(x))dx\|\leq\frac{C_k}{n}
\end{equation}
for large enough natural $n$.
Indeed, let $p:\tor\rightarrow SU(2)$ be a measurable function such that
\[p(x)\varphi(x)p(Tx)^{-1}=\delta(x)=
\left[\begin{array}{cc}
\gamma(x) & 0\\
0 & \overline{\gamma(x)}
\end{array}\right].\]
Consider the unitary operator
$V:\Hh_1^{\Pi_{2k-1}}\rightarrow\Hh_1^{\Pi_{2k-1}}$ given by
\[V(\sum_{i=1}^{d_{\Pi_{2k-1}}}\Pi_{1i}(g)f_i(x))=
\sum_{i,j=1}^{d_{\Pi_{2k-1}}}\Pi_{1j}(g)\Pi_{ji}(p(x)^{-1})f_i(x)).\]
Then
\begin{eqnarray*}
V^{-1}U_{T_{\varphi}}V(\sum_{i=1}^{d_{\Pi_{2k-1}}}\Pi_{1i}(g)f_i(x)) & = &
\sum_{i,j,l,m=1}^{d_{\Pi_{2k-1}}}\Pi_{1m}(g)\Pi_{ml}(p(x))
\Pi_{lj}(\varphi(x))\Pi_{ji}(p(Tx)^{-1})f_i(Tx)\\
& = &
\sum_{i=1}^{d_{\Pi_{2k-1}}}\Pi_{1i}(g)
\Pi_{ii}(\delta(x))f_i(Tx).
\end{eqnarray*}
From (\ref{lebe}),
$U_{T_{\varphi}}:\Hh_1^{\Pi_{2k-1}}\rightarrow\Hh_1^{\Pi_{2k-1}}$ 
has Lebesgue spectrum of uniform multiplicity, by Corollary~\ref{kica}.
Hence $V^{-1}U_{T_{\varphi}}V$ 
has Lebesgue spectrum of uniform multiplicity and it is the product of
the operators $U_j:L^2(\tor,{\Bbb C})\rightarrow L^2(\tor,{\Bbb C})$ given
by $U_jf(x)=(\gamma(x))^{2k-2j+1}f(Tx)$ for $j=1,\ldots,2k$. Therefore $U_j$
has absolutely continuous spectrum for $j=1,\ldots,2k$.
By Lemma~\ref{unif}, $U_j$ has Lebesgue spectrum for all $j=1,\ldots,2k$
and $k\in{\Bbb N}$.
It follows that the Lebesgue component in the spectrum of
$T_{\gamma(\varphi)}$ has countable multiplicity.

By assumption, 
\[\|\frac{1}{n}D\varphi^{(n)}(\varphi^{(n)})^{-1}\|\rightarrow d(\varphi)\]
uniformly. Therefore
\begin{equation}\label{cik}
\|\frac{1}{n}D\varphi^{(n)}(x)(\varphi^{(n)}(x))^{-1}\|\geq d(\varphi)/2
\end{equation}
for large enough natural $n$.
For all $A,B\in M_k({\Bbb C})$ we have $\|AB\|\leq\sqrt{k}\|A\|\|B\|$. Applying
these facts, (\ref{cal1}) and (\ref{cal2}) we get 
\begin{eqnarray*}
\|\int_{\tor}\Pi_{2k-1}(\varphi^{(n)}(x))dx\| & = & 
\|\int_{\tor}\Pi_{2k-1}(\varphi^{(n)}(x))
(D\Pi_{2k-1}(\varphi^{(n)}(x)))^{-1}d\Pi_{2k-1}(\varphi^{(n)}(x))\| \\
 & = &
\|\int_{\tor}(\Pi_{2k-1}^*(D\varphi^{(n)}(x)(\varphi^{(n)}(x))^{-1}))^{-1}
d\Pi_{2k-1}(\varphi^{(n)}(x))\| \\
 & = &
\|\int_{\tor}(\Pi_{2k-1}(\varphi^{(n)}(x)))^T
d((\Pi_{2k-1}^*(D\varphi^{(n)}(x)(\varphi^{(n)}(x))^{-1}))^{-1})^T\| \\
 & = &
\|\int_{\tor}[(\Pi_{2k-1}(\varphi^{(n)}(x)))^T
((\Pi_{2k-1}^*(D\varphi^{(n)}(x)(\varphi^{(n)}(x))^{-1}))^{-1})^T\\
 & & 
\mbox{\hspace{7mm}}
(\Pi_{2k-1}^*D(D\varphi^{(n)}(x)(\varphi^{(n)}(x))^{-1}))^T\\
 & & 
\mbox{\hspace{7mm}}
((\Pi_{2k-1}^*(D\varphi^{(n)}(x)(\varphi^{(n)}(x))^{-1}))^{-1})^T]dx\| \\
 & \leq &
2k\int_{\tor}
[\|(\Pi_{2k-1}^*(D\varphi^{(n)}(x)(\varphi^{(n)}(x))^{-1}))^{-1}\|^2\\
 & & 
\mbox{\hspace{10mm}}
\|\Pi_{2k-1}^*D(D\varphi^{(n)}(x)(\varphi^{(n)}(x))^{-1})\|]dx.
\end{eqnarray*}
By Lemma~\ref{kickic}, we have
\[\|(\Pi_{2k-1}^*(D\varphi^{(n)}(x)(\varphi^{(n)}(x))^{-1}))^{-1}\|\leq
K_k\|D\varphi^{(n)}(x)(\varphi^{(n)}(x))^{-1}\|^{-1}.\]
From this, (\ref{taki}) and (\ref{cik}) we obtain
\begin{eqnarray*}
\lefteqn{\|\int_{\tor}\Pi_{2k-1}(\varphi^{(n)}(x))dx\|}\\
& \leq &
\frac{K_k^2(2k)^3}{n}\int_{\tor}
[\|\frac{1}{n}D\varphi^{(n)}(x)(\varphi^{(n)}(x))^{-1})\|^{-2}
\|D(\frac{1}{n}D\varphi^{(n)}(x)(\varphi^{(n)}(x))^{-1})\|]dx\\
& \leq & 
\frac{1}{n}(\frac{8K_kk^2}{d(\varphi)})^2
\|D(\frac{1}{n}D\varphi^{(n)}(\varphi^{(n)})^{-1})\|_{L^1}
\end{eqnarray*}
for large enough natural $n$. By assumption, there exists a real constant
$M>0$ such that
$\|D(\frac{1}{n}D\varphi^{(n)}(\varphi^{(n)})^{-1})\|_{L^1}\leq M$. 
Then
\[\|\int_{\tor}\Pi_{2k-1}(\varphi^{(n)}(x))dx\|\leq\frac{C_k}{n}\]
for large enough natural $n$, where $C_k=(\frac{8K_kk^2}{d(\varphi)})^2M$.
$\blacksquare$ 

\vspace{2ex}

In this section we also present a class of cocycles satisfying the
assumptions of Theorem~\ref{glowne1}. We will need the following lemma.

\begin{lem}
Let $\{f_n:\tor\rightarrow{\Bbb C}^d;n\in{\Bbb N}\}$ be a sequence of
absolutely continuous functions. Assume that the
sequence $\{f_n\}_{n\in{\Bbb N}}$ converges in $L^1(\tor,{\Bbb R}^d)$ to a
function $f$ and it is  bounded
for the sup norm.  Suppose that there
is a sequence $\{h_n\}_{n\in{\Bbb N}}$ convergent in $L^2_+(\tor,{\Bbb R})$
and a sequence $\{k_n\}_{n\in{\Bbb N}}$ bounded in $L^2_+(\tor,{\Bbb R})$
such that
\[\|Df_n(x)\|\leq h_n(x)k_n(x)\;\;\mbox{for $\lambda$--a.e.\ $x\in\tor$}\]
and for any natural $n$. Then $\{f_n\}_{n\in{\Bbb N}}$ converges to
$f$ uniformly. 
\end{lem}

\pf
Denote by $h\in L^2_+(\tor,{\Bbb R})$ the limit of the  sequence
$\{h_n\}_{n\in{\Bbb N}}$. Let $M>0$ be a real number such that
$\|k_n\|_{L^2}\leq M$ for all natural $n$. First observe that the sequence
$\{f_n\}_{n\in{\Bbb N}}$ is equicontinuous. Fix $\ep>0$. Take
$n_0\in{\Bbb N}$ such that $\|h_n-h\|_{L^2}<\ep/2M$ for any $n\geq n_0$.
Then for all $x,y\in\tor$ and $n\geq n_0$ we have
\begin{eqnarray*}
\|f_n(x)-f_n(y)\| & = & \|\int_x^yDf_n(t)dt\|\leq \int_x^y\|Df_n(t)\|dt\\
 & \leq & \int_x^yh_n(t)k_n(t)dt \leq \|k_n\|_{L^2}\sqrt{\int_x^yh_n^2(t)dt}\\
 & \leq & M(\sqrt{\int_x^yh^2(t)dt}+\|h_n-h\|_{L^2})\leq 
M(\sqrt{\int_x^yh^2(t)dt}+\frac{\ep}{2M}).
\end{eqnarray*}
Choose $\delta_1>0$ such that $|x-y|<\delta_1$ implies
$\int_x^yh^2(t)dt<(\ep/2M)^2$. Hence if $|x-y|<\delta_1$, then
$\|f_n(x)-f_n(y)\|<\ep$ for any $n\geq n_0$. Next choose
$0<\delta\leq\delta_1$ such that $|x-y|<\delta$ implies
$\|f_n(x)-f_n(y)\|<\ep$ for any $n\leq n_0$. It follows that if
$|x-y|<\delta$, then $\|f_n(x)-f_n(y)\|<\ep$ for every natural $n$.

By the Arzela--Ascoli theorem, for any subsequence of $\{f_n\}_{n\in{\Bbb N}}$ there
exists a subsequence convergent to $f$ uniformly. Consequently, the 
sequence $\{f_n\}_{n\in{\Bbb N}}$ converges to $f$ uniformly. $\blacksquare$

\vspace{2ex}

This gives the following corollary.

\begin{cor}\label{wnios}
Let $\{f_n:\tor\rightarrow{\Bbb C}^d;n\in{\Bbb N}\}$ be a sequence of
absolutely continuous functions. Assume that the
sequence $\{f_n\}_{n\in{\Bbb N}}$ converges in $L^1(\tor,{\Bbb R}^d)$ to a
function $f$ and it is bounded
for the sup norm. Suppose that there
is a sequence $\{h_n\}_{n\in{\Bbb N}}$ convergent in $L^2_+(\tor,{\Bbb R})$,
a sequence $\{k_n\}_{n\in{\Bbb N}}$ bounded in $L^2_+(\tor,{\Bbb R})$ and
a sequence $\{l_n\}_{n\in{\Bbb N}}$ convergent in $L^1_+(\tor,{\Bbb R})$
such that
\[\|Df_n(x)\|\leq l_n(x)+ h_n(x)k_n(x)\;\;\mbox{for $\lambda$-a.e.
$x\in\tor$}\] 
and for any natural $n$. Then $\{f_n\}_{n\in{\Bbb N}}$ converges to
$f$ uniformly. $\blacksquare$
\end{cor}

We will denote by $BV^2(\tor,SU(2))$ the set of all functions
$f:\tor\rightarrow SU(2)$ of bounded variation 
such that $Df(f)^{-1}\in L^2(\tor,\su)$.

\begin{lem}\label{osta}
Let $\varphi:\tor\rightarrow SU(2)$ be a $C^2$--cocycle. Suppose
that $\varphi$ is cohomologous to a diagonal cocycle with a transfer
function in $BV^2(\tor,SU(2))$.
Then the sequence
$\{\frac{1}{n}D\varphi^{(n)}(\varphi^{(n)})^{-1})\}_{n\in{\Bbb N}}$
is uniformly convergent and
$\{D(\frac{1}{n}D\varphi^{(n)}(\varphi^{(n)})^{-1})\}_{n\in{\Bbb N}}$
is bounded in $L^1(\tor,\su)$ 
\end{lem}

\pf 
By Corollary~\ref{wnios}, it suffices to show that there exist
a sequence $\{h_n\}_{n\in{\Bbb N}}$ convergent in $L^2_+(\tor,{\Bbb R})$,
a sequence $\{k_n\}_{n\in{\Bbb N}}$ bounded in $L^2_+(\tor,{\Bbb R})$ and
a sequence $\{l_n\}_{n\in{\Bbb N}}$ convergent in $L^1_+(\tor,{\Bbb R})$
such that
\[\|D(\frac{1}{n}D\varphi^{(n)}(x)(\varphi^{(n)}(x))^{-1})\|
\leq l_n(x)+ h_n(x)k_n(x)\;\;\mbox{for $\lambda$-a.e.\ $x\in\tor$}.\]
By assumption, there exist $\delta,p\in BV^2(\tor,SU(2))$ such that
$\varphi(x)=p(x)^{-1}\delta(x)p(Tx)$, where $\delta$ is a diagonal cocycle.
Then
\[D\varphi(x)\varphi(x)^{-1}=
-p(x)^{-1}Dp(x)+p(x)^{-1}D\delta(x)\delta(x)^{-1}p(x)+
\varphi(x)p(Tx)^{-1}Dp(Tx)\varphi(x)^{-1}\]
for $\lambda$--a.e.\ $x\in\tor$. Set 
\[\tv(x)=D\varphi(x)\varphi(x)^{-1},\;\;\tp(x)=p(x)^{-1}Dp(x)\;\mbox{ and }
\;\td(x)=p(x)^{-1}D\delta(x)\delta(x)^{-1}p(x).\] 
Then $\tv(x)=-\tp(x)+U\tp(x)+\td(x)$, where $\tp,\td\in L^2(\tor,\su)$.
Since
\[\frac{1}{n}D\varphi^{(n)}(\varphi^{(n)})^{-1}=
\frac{1}{n}\sum_{k=0}^{n-1}\varphi^{(k)}\,\tv\circ T^k\,(\varphi^{(k)})^{-1},\]
we have
\begin{eqnarray*}
\lefteqn{
D(\frac{1}{n}D\varphi^{(n)}(\varphi^{(n)})^{-1})}\\ 
& = &
\frac{1}{n}\sum_{k=0}^{n-1}\sum_{j=0}^{k-1}
(\ad_{\varphi^{(j)}}\,(\tv\circ T^j)\ad_{\varphi^{(k)}}
\,(\tv\circ T^k)
-\ad_{\varphi^{(k)}}\,(\tv\circ T^k)\ad_{\varphi^{(j)}}\,(\tv\circ T^j))\\
& & + \frac{1}{n}\sum_{k=0}^{n-1}\ad_{\varphi^{(k)}}\,(D\tv\circ T^k)\\
 & = & \frac{1}{n}\sum_{k=0}^{n-1}\sum_{j=0}^{k-1}[U^j\tv,U^k\tv]+
\frac{1}{n}\sum_{k=0}^{n-1}U^k(D\tv).
\end{eqnarray*}
However,
\begin{eqnarray*}
\sum_{k=0}^{n-1}\sum_{j=0}^{k-1}[U^j\tv,U^k\tv] & = &
\sum_{k=0}^{n-1}\sum_{j=0}^{k-1}[U^{j+1}\tp-U^j\tp+U^j\td,U^k\tv]\\
 & = & 
\sum_{k=0}^{n-1}[U^{k}\tp-\tp,U^k\tv]+
\sum_{k=0}^{n-1}\sum_{j=0}^{k-1}[U^j\td,U^{k+1}\tp-U^k\tp+U^k\td]\\
 & = & 
\sum_{k=0}^{n-1}[U^{k}\tp-\tp,U^k\tv]+
\sum_{k=0}^{n-1}\sum_{j=0}^{k-1}[U^j\td,U^k\td]+
\sum_{j=0}^{n-2}[U^j\td,U^{n}\tp-U^{j+1}\tp].
\end{eqnarray*}
Since
\begin{eqnarray*}
U^j\td(x) & = &
\ad_{\varphi^{(n)}(x)p(T^nx)^{-1}}(D\delta(T^nx)\delta(T^nx)^{-1})\\
 & = &
\ad_{p(x)^{-1}\delta^{(n)}(x)}(D\delta(T^nx)\delta(T^nx)^{-1})\\
 & = &
\ad_{p(x)^{-1}}(D\delta(T^nx)\delta(T^nx)^{-1}),
\end{eqnarray*}
we have $[U^j\td,U^k\td]=0$ for any integers $j,k$. 
Observe that $\|[A,B]\|\leq 2\|A\|\|B\|$  for any $A,B\in\su$.
It follows that
\begin{eqnarray*}
\lefteqn{
\|D(\frac{1}{n}D\varphi^{(n)}(\varphi^{(n)})^{-1})\|}\\ 
& \leq &
\frac{2}{n}\sum_{k=0}^{n-1}(\|D\tv\circ T^k\|+\|\tp\circ T^k\|
\|\tv\circ T^k\|
+\|\tp\|\|\tv\circ T^k\|+\|\td\circ T^k\|\|\tp\circ T^{k+1}\|)\\
 & & +
\|\tp\circ T^n\|\frac{2}{n}\sum_{k=0}^{n-1}\|\td\circ T^k\|.
\end{eqnarray*}
Set
\begin{eqnarray*}
h_n & = &
\frac{2}{n}\sum_{k=0}^{n-1}\|\td\circ T^k\|\\
k_n & = & \|\tp\circ T^n\|\\
l_n & = & 
\frac{2}{n}\sum_{k=0}^{n-1}(\|D\tv\circ T^k\|+\|\tp\circ T^k\|
\|\tv\circ T^k\|
+\|\tp\|\|\tv\circ T^k\|+\|\td\circ T^k\|\|\tp\circ T^{k+1}\|). 
\end{eqnarray*}
By the Birkhoff ergodic
theorem, the sequence $\{h_n\}_{n\in{\Bbb N}}$ converges in
$L^2_+(\tor,{\Bbb R})$ and the
sequence $\{l_n\}_{n\in{\Bbb N}}$ converges in $L^1_+(\tor,{\Bbb R})$.
This completes the proof. $\blacksquare$

\vspace{2ex}

Theorem~\ref{glowne1} and Lemma~\ref{osta} lead to the following conclusion.

\begin{cor}\label{glowne2}
Let $\varphi:\tor\rightarrow SU(2)$ be a $C^2$--cocycle  with
$d(\varphi)\neq 0$. Suppose that 
$\varphi$ is cohomologous to a diagonal cocycle with a transfer
function in $BV^2(\tor,SU(2))$.
Then  the Lebesgue component in the spectrum of
$T_{\varphi}$ has countable multiplicity. Moreover, the Lebesgue component
in the spectrum of $T_{\gamma(\varphi)}$ has countable multiplicity, too.
\end{cor}

The following result will be useful in the next section of the paper.

\begin{prop}\label{prop}
For every $C^2$--cocycle $\varphi:\tor\rightarrow SU(2)$, the sequence
\[\frac{1}{n^2}D(D\varphi^{(n)}(\varphi^{(n)})^{-1})\] 
converges to zero in $L^1(\tor,\su)$.
\end{prop}

The following lemmas are some simple generalizations of some classical results.
Their proofs are left to the reader.

\begin{lem}\label{stoltz}
Let $\{a_n\}_{n\in{\Bbb N}}$ be an increasing sequence of natural numbers and
let $\{f_n\}_{n\in{\Bbb N}}$ be a sequence in the Banach space
$L^2(\tor, M_2({\Bbb C}))$. Then 
\[\frac{f_{n+1}-f_n}{a_{n+1}-a_n}\rightarrow g\mbox{ in }
L^2(\tor, M_2({\Bbb C}))\;\;\Longrightarrow\;\;\frac{f_n}{a_n}\rightarrow g\mbox{ in }
L^2(\tor, M_2({\Bbb C})).\]
\end{lem}

\begin{lem}\label{prop2}
Let $\{g_k^n;n\in{\Bbb N}, 0\leq k<n\}$ be a triangular matrix of elements
from $L^2(\tor, M_2({\Bbb C}))$ such that $\|g_k^n\|=O(1/n)$ and
\[g_0^n+g_1^n+\ldots+g_{n-1}^n\rightarrow g \;\;\mbox{ in }\;\;L^2(\tor, M_2({\Bbb C})).\]
Then $f_n\rightarrow f$ in $L^2(\tor, M_2({\Bbb C}))$ implies
\[\sum_{k=0}^{n-1}g_k^n\,f_k\rightarrow g\,f\;\;\mbox{ and }\;\;
\sum_{k=0}^{n-1}f_k\,g_k^n\rightarrow f\,g\;\;\mbox{ in }\;\;L^1(\tor, M_2({\Bbb C})).\]
\end{lem} 

{\bf Proof of Proposition~\ref{prop}.}
First recall that
\[\frac{1}{n^2}D(D\varphi^{(n)}(\varphi^{(n)})^{-1})
=\frac{1}{n^2}\sum_{k=0}^{n-1}\sum_{j=0}^{k-1}[U^j\tv,U^k\tv]+
\frac{1}{n^2}\sum_{k=0}^{n-1}U^k(D\tv),\]
where $\tv=D\varphi(\varphi)^{-1}$ and
\[\frac{1}{n}\sum_{k=0}^{n-1}U^k\tv\rightarrow \psi\;\;\mbox{ in }
L^2(\tor,\su).\]
Since $\frac{1}{n^2}\sum_{k=0}^{n-1}U^k(D\tv)$ uniformly
converges to zero, it suffices to show that
\[\lim_{n\rightarrow\infty}\frac{1}{n^2}
\sum_{k=0}^{n-1}\sum_{j=0}^{k-1}U^j\tv\,U^k\tv=
\lim_{n\rightarrow\infty}\frac{1}{n^2}
\sum_{k=0}^{n-1}\sum_{j=0}^{k-1}U^k\tv\,U^j\tv=\frac{1}{2}\psi\,\psi
\;\;\mbox{ in }\;\;
L^2(\tor, M_2({\Bbb C})).\]
Set $f_n=\sum_{k=0}^{n-1}(n-k)U^k\tv$ and $a_n=n^2$. Then
\[\frac{f_{n+1}-f_n}{a_{n+1}-a_n}=
\frac{\sum_{k=0}^{n}(n+1-k)U^k\tv-\sum_{k=0}^{n-1}(n-k)U^k\tv}{(n+1)^2-n^2}=
\frac{\sum_{k=0}^{n}U^k\tv}{2n+1}\rightarrow \frac{1}{2}\psi\]
in $L^2(\tor, M_2({\Bbb C}))$. 
Applying Lemma~\ref{stoltz},  we get
\[\frac{1}{n^2}\sum_{k=0}^{n-1}(n-k)U^k\tv\rightarrow \frac{1}{2}
\psi\;\;\mbox{ in }\;\;L^2(\tor, M_2({\Bbb C})).\]
Therefore
\[\frac{1}{n^2}\sum_{k=0}^{n-1}k\,U^k\tv=
\frac{1}{n}\sum_{k=0}^{n-1}U^k\tv-
\frac{1}{n^2}\sum_{k=0}^{n-1}(n-k)U^k\tv\rightarrow\psi-\frac{1}{2}\psi=
\frac{1}{2}\psi\]
in $L^2(\tor, M_2({\Bbb C}))$.
Applying Lemma~\ref{prop2} with $g_k^n=\frac{k}{n^2}U^k\tv$ and 
$f_k=\frac{1}{k}\sum_{j=0}^{k-1}U^j\tv$, we conclude that
\[\sum_{k=0}^{n-1}g_k^n\,f_k=
\frac{1}{n^2}\sum_{k=0}^{n-1}\sum_{j=0}^{k-1}U^k\tv\,U^j\tv
\rightarrow \frac{1}{2}\psi\,\psi\]
and
\[\sum_{k=0}^{n-1}f_k\,g_k^n=\
\frac{1}{n^2}\sum_{k=0}^{n-1}\sum_{j=0}^{k-1}U^j\tv\,U^k\tv
\rightarrow \frac{1}{2}\psi\,\psi\]
in $L^2(\tor, M_2({\Bbb C}))$, which completes the proof. $\blacksquare$

\section{Possible values of degree}\indent

\label{integ}
One may ask what we know about the set of possible values of degree. 
For $G=\tor$ the degree of each 
smooth cocycle is an integer number. Probably, in the case of cocycles
with values in $SU(2)$ the set of 
possible values of degree is more complicated. However, in this
section we show that if $\alpha$ is the golden ratio, then
the degree of each smooth cocycle  belongs to $2\pi{\Bbb N}_0$. The idea
of renormalization,
which is used to prove this result is due to Rychlik \cite{Ry}.

Let $\alpha$ be the golden ratio (i.e.\ the positive root of the equation
$\alpha^2+\alpha=1$). It will be advantageous for our notation to consider
the interval $[-\alpha^2,\alpha)$ to be the model of the circle. Then  the
map $T:[-\alpha^2,\alpha)\rightarrow[-\alpha^2,\alpha)$ given by
\[T(x)=\left\{
\begin{array}{lcl}
x+\alpha & \mbox{ for } & x\in[-\alpha^2,0)\\
x-\alpha^2 & \mbox{ for } & x\in[0,\alpha)
\end{array}\right.\]
is the rotation by $\alpha$. Let $X=[-\alpha^2,\alpha^3)$. Then the first
return time to $X$, which we call $\tau$, satisfies the following formula
\[\tau(x)=\left\{
\begin{array}{rcl}
1 & \mbox{ for } & x\in[0,\alpha^3)\\
2 & \mbox{ for } & x\in[-\alpha^2,0)
\end{array}\right.\]
and the first return map $T_X:X\rightarrow X$ is equal to $T$ up to a
linear scaling. Indeed, if $M:\tor\rightarrow X$ is the map given by
$M(x)=-\alpha x$, then $T_X\circ M=M\circ T$.\\
By $W^1$ we mean the space of all cocycle $\varphi:\tor\rightarrow SU(2)$
such that the functions $\varphi:[-\alpha^2,0)\rightarrow SU(2)$, 
$\varphi:[0,\alpha)\rightarrow SU(2)$ are both of class $C^1$ 
and 
\[\lim_{x\rightarrow 0^-}D\varphi(x)\varphi(x)^{-1}\;\;\mbox{ and }\;\;
\lim_{x\rightarrow \alpha^-}D\varphi(x)\varphi(x)^{-1} \]
exist. The topology of $W^1$ is induced from
$C^1((-\alpha^2,0)\cup(0,\alpha))$. 
Consider the renormalization operator $\Phi:W^1\rightarrow W^1$ defined by
\[\Phi\:\varphi(x)=\varphi^{(\tau(Mx))}(Mx).\]
Then
\[\Phi^n \varphi(x) =\left\{
\begin{array}{rcl}
\varphi^{(q_{n+1})} (M^nx) & \mbox{ for } & x\in[-\alpha^2,0)\\
\varphi^{(q_{n+2})} (M^nx) & \mbox{ for } & x\in[0,\alpha)
\end{array}\right.\]
for any natural $n$, where $\{q_n\}_{n\in{\Bbb N}}$ is the Fibonacci
sequence. 
By $W_0^1$ we mean the set of all
cocycles $\varphi\in W^1$ such that $\varphi^{(2)}$ is continuous at $0$.
The set $W^1_0$ is a closed subset of $W^1$ and
\begin{equation}
\Phi(W^1_0)\subset W^1_0
\end{equation} 
(see \cite{Ry}). 
It is easy to check that
$\|D(\Phi\varphi)(\Phi\varphi)^{-1}\|_{L^1}\leq
\|D\varphi(\varphi)^{-1}\|_{L^1}$ for any $\varphi\in W^1$.   
The following result is due to M.\ Rychlik~\cite{Ry}.

\begin{prop}\label{rych}
If $\|D(\Phi^k\varphi)(\Phi^k\varphi)^{-1}\|_{L^1}=
\|D\varphi(\varphi)^{-1}\|_{L^1}$ for all natural $k$, then
\[D\varphi(x)(\varphi(x))^{-1}=\alpha\,\mbox{\em Ad}_{\varphi(x)}
[D\varphi(Tx)(\varphi(Tx))^{-1}]\]
for every $x\in [-\alpha^2,0)$. $\blacksquare$
\end{prop}

\begin{lem}\label{pomost}
Let $\varphi:\tor\rightarrow SU(2)$ be a $C^2$--cocycle. Assume that 
\[\frac{1}{n}D\varphi^{(n)}(0)(\varphi^{(n)}(0))^{-1}\rightarrow H\in\su\]
and there is an increasing sequence $\{n_k\}_{k\in{\Bbb N}}$ of even
numbers such that 
\[\lim_{k\rightarrow\infty}\alpha^{n_k}\int_0^{\alpha^{n_k}}
|D(D\varphi^{(q_{n_k+i})}(x)(\varphi^{(q_{n_k+i})}(x))^{-1})|dx=0\]
for $i=1,2$. Then
$\|H\|\in 2\pi{\Bbb N}_0$.
\end{lem}

\pf
First note that
\[D\Phi^n \varphi(x)(\Phi^n \varphi(x))^{-1} =\left\{
\begin{array}{rcl}
\alpha^n D\varphi^{(q_{n+1})} (M^nx)
(\varphi^{(q_{n+1})} (M^nx))^{-1} & \mbox{ for } & x\in[-\alpha^2,0)\\
\alpha^n D\varphi^{(q_{n+2})} (M^nx)
(\varphi^{(q_{n+2})} (M^nx))^{-1} & \mbox{ for } & x\in[0,\alpha)
\end{array}\right.\]
for any even $n$.
Since 
\begin{eqnarray*}
\lefteqn{
|\frac{1}{q_{n+i}}D\varphi^{(q_{n+i})}(M^nx)
(\varphi^{(q_{n+i})}(M^nx))^{-1}-
\frac{1}{q_{n+i}}D\varphi^{(q_{n+i})}(0)
(\varphi^{(q_{n+i})}(0))^{-1}|}
\\
 & \leq &\frac{1}{q_{n+i}}\int_0^{\alpha^nx}
|D(D\varphi^{(q_{n+i})}(\varphi^{(q_{n+i})})^{-1})|d\lambda
  \leq  \frac{1}{q_{n+i}\alpha^n}\alpha^n\int_0^{\alpha^n}
|D(D\varphi^{(q_{n+i})}(\varphi^{(q_{n+i})})^{-1})|d\lambda
\end{eqnarray*} 
for all even $n$, $i=1,2$ and 
\[\lim_{n\rightarrow\infty}\alpha^n q_{n+1}=1/(1+\alpha^2),\;\;\;
\lim_{n\rightarrow\infty}\alpha^n q_{n+2}=1/(\alpha+\alpha^3),\]
it follows that
\begin{eqnarray*}
\lim_{k\rightarrow\infty}
D\Phi^{n_k} \varphi(x)(\Phi^{n_k} \varphi(x))^{-1} & = &
\lim_{k\rightarrow\infty}
\alpha^{n_k}q_{n_k+1}\frac{1}{q_{n_k+1}}D\varphi^{(q_{n_k+1})}(0)
(\varphi^{(q_{n_k+1})}(0))^{-1}\\
 & = & \frac{1}{1+\alpha^2}H
\end{eqnarray*}
uniformly on $[-\alpha^2,0)$ and
\begin{eqnarray*}
\lim_{k\rightarrow\infty}
D\Phi^{n_k} \varphi(x)(\Phi^{n_k} \varphi(x))^{-1} & = & 
\lim_{k\rightarrow\infty}
\alpha^{n_k}q_{n_k+2}\frac{1}{q_{n_k+2}}D\varphi^{(q_{n_k+2})}(0)
(\varphi^{(q_{n_k+2})}(0))^{-1}\\
 & = & \frac{1}{\alpha+\alpha^3}H
\end{eqnarray*}
uniformly on $[0,\alpha)$.
Therefore we can assume that there exists $v\in W^1$ such that
\[\Phi^{n_k}\varphi\rightarrow v\;\;\mbox{ and }\;\;
D\Phi^{n_k} \varphi(\Phi^{n_k} \varphi)^{-1}\rightarrow Dv\:v^{-1}\]
uniformly.  Then
\[Dv(x)(v(x))^{-1} =\left\{
\begin{array}{rcl}
\alpha A & \mbox{ for } & x\in[-\alpha^2,0)\\
A & \mbox{ for } & x\in[0,\alpha),
\end{array}\right.\]
where  $A=1/(\alpha+\alpha^3)\,H\in\su$.
Therefore 
\[v(x)=\left\{
\begin{array}{rcl}
e^{\alpha xA}B & \mbox{ for } & x\in[-\alpha^2,0)\\
e^{xA}C & \mbox{ for } & x\in[0,\alpha),
\end{array}\right.\]
where $B=v_-(0)$ and $C=v_+(0)$. 
Since the set $W^1_0\subset W^1$ is closed and $\Phi$--invariant, $v\in
W^1_0$.  It follows that 
\begin{equation}\label{pich}
Ce^{-\alpha^3A}B=Be^{\alpha A}C.
\end{equation}
Since $v$ is a limit point of the sequence
$\{\Phi^n\varphi\}_{n\in{\Bbb N}}$, we have 
$\|D\Phi^kv(\Phi^kv)^{-1}\|_{L^1}=
\|Dvv^{-1}\|_{L^1}$ for any natural $k$. By Proposition~\ref{rych},
\[\lim_{x\rightarrow 0^-}Dv(x)(v(x))^{-1}=\alpha\mbox{ Ad}_{v_-(0)}
\lim_{x\rightarrow \alpha^-}Dv(x)(v(x))^{-1}.\]
Hence
\[\alpha A=\alpha\,\mbox{Ad}_{B}(A)\]
and finally $AB=BA$. 
Therefore
\[\Phi v(x)=\left\{
\begin{array}{rcl}
e^{-\alpha xA}C & \mbox{ for } & x\in[-\alpha^2,0)\\
e^{-xA+\alpha A}BC & \mbox{ for } & x\in[0,\alpha).
\end{array}\right.\]
By Proposition~\ref{rych},
\[\lim_{x\rightarrow 0^-}D\Phi v(x)(\Phi v(x))^{-1}=\alpha\mbox{ Ad}_{\Phi v_-(0)}
\lim_{x\rightarrow \alpha^-}D\Phi v(x)(\Phi v(x))^{-1}.\]
Hence
\[-\alpha\, A=\alpha\,\mbox{Ad}_{C}(-A)\]
and finally $AC=CA$.
It follows that $B$ and $C$ commute, by (\ref{pich}).
From~(\ref{pich}), we obtain $e^{(\alpha+\alpha^3)A}=\mbox{Id}$. 
Therefore
$\|H\|=\|(\alpha+\alpha^3)A\|\in 2\pi{\Bbb N}_0$. $\blacksquare$ 

\begin{theo}
Suppose that $\alpha$ is the golden ratio. Then for every $C^2$--cocycle
$\varphi:\tor\rightarrow SU(2)$, we have $d(\varphi)\in 2\pi{\Bbb N}_0$.  
\end{theo}

\pf
Fix $n\in{\Bbb N}$ such that $2\alpha^{2n}[1/2\alpha^{2n}]\geq 4/5$. Set
$I_j=[2(j-1)\alpha^{2n},2j\alpha^{2n}]$ for
$j\in E=\{1,\ldots,[1/2\alpha^{2n}]\}$ and $\ep_n=\frac{1}{q_n^2}
\int_{\tor}|D(D\varphi^{(q_n)}(\varphi^{(q_n)})^{-1})|d\lambda$.
By Proposition~\ref{prop}, $\ep_n$ tends to zero.
For $i=1,2$ define
\[E_i=\{j\in E;\;
\frac{1}{2\alpha^{2n}q_{2n+i}^2}
\int_{I_j}|D(D\varphi^{(q_{2n+i})}(\varphi^{(q_{2n+i})})^{-1})|d\lambda
\leq 10\ep_{2n+i}\}.\]
Then
\begin{eqnarray*}
\ep_{2n+i} & = & \frac{1}{q_{2n+i}^2}
\int_{\tor}|D(D\varphi^{(q_{2n+i})}(\varphi^{(q_{2n+i})})^{-1})|d\lambda\\
& \geq & \frac{1}{q_{2n+i}^2}
\sum_{j\in E\setminus E_i}\int_{I_j}
|D(D\varphi^{(q_{2n+i})}(\varphi^{(q_{2n+i})})^{-1})|d\lambda\\
& \geq & 20\alpha^{2n}\ep_{2n+i}([1/2\alpha^{2n}]-\#E_i).
\end{eqnarray*}
Hence
\[\#E_i\geq[1/2\alpha^{2n}]
(1-\frac{1}{10}\frac{1/2\alpha^{2n}}{[1/2\alpha^{2n}]})\geq 
\frac{7}{8}[1/2\alpha^{2n}]\]
for $i=1,2$. Therefore
\[\#(E_1\cap E_2)\geq\#E_1+\#E_2-\#E\geq \frac{3}{4}[1/2\alpha^{2n}].\]
Define
\[G_n=\bigcup_{j\in E_1\cap E_2}[(2j-2)\alpha^{2n},(2j-1)\alpha^{2n}].\]
Observe that  $y\in G_n$ implies
\[\frac{1}{2\alpha^{2n}q_{2n+i}^2}
\int_{y}^{y+\alpha^{2n}}
|D(D\varphi^{(q_{2n+i})}(\varphi^{(q_{2n+i})})^{-1})|d\lambda
\leq 10\ep_{2n+i}\]
for $i=1,2$ and
\[\lambda(G_n)\geq\alpha^{2n}\#(E_1\cap E_2)\geq
\frac{3}{8}2\alpha^{2n}[1/2\alpha^{2n}]\geq\frac{3}{10}.\]
Set $\displaystyle G'=\bigcap_{n\in{\Bbb N}}\bigcup_{k>n}G_k$. Then
$\lambda(G')\geq 3/10$. Since
$\frac{1}{n}D\varphi^{(n)}(\varphi^{(n)})^{-1}\rightarrow \psi$ almost
everywhere, we see that the set
\[G=\{x\in G';\frac{1}{n}D\varphi^{(n)}(x)(\varphi^{(n)}(x))^{-1}
\rightarrow \psi(x)\}\]
has positive measure.

For every $y\in\tor$ denote by $\varphi_y:\tor\rightarrow SU(2)$ the
$C^2$--cocycle $\varphi_y(x)=\varphi(x+y)$.
Suppose that $y\in G$. Then $\frac{1}{n}D\varphi_y^{(n)}(0)
(\varphi_y^{(n)}(0))^{-1}\rightarrow \psi(y)$ and there exists an increasing
sequence $\{n_k\}_{k\in{\Bbb N}}$ of natural numbers such that $y\in
G_{n_k}$ for any natural $k$. Hence
\[\alpha^{2n_k}\int_{0}^{\alpha^{2n_k}}
|D(D\varphi_y^{(q_{2n_k+i})}(\varphi_y^{(q_{2n_k+i})})^{-1})|d\lambda
\leq 20(\alpha^{2n_k}q_{2n_k+i})^2\ep_{2n_k+i}\]
for $i=1,2$. Since the sequence $\{\alpha^nq_{n+i}\}_{n\in{\Bbb N}}$
converges for $i=1,2$ and $\ep_n$ tends to zero, letting
$k\rightarrow\infty$ we have 
\[\lim_{k\rightarrow\infty}\alpha^{2n_k}\int_{0}^{\alpha^{2n_k}}
|D(D\varphi_y^{(q_{2n_k+i})}(\varphi_y^{(q_{2n_k+i})})^{-1})|d\lambda=0\]
for $i=1,2$.
By Lemma~\ref{pomost}, $\|\psi(y)\|\in 2\pi {\Bbb N}_0$ for every $y\in G$.
Since $d(\varphi)=\|\psi(y)\|$ for a.e.\ $y\in\tor$, we conclude that
$d(\varphi)\in 2\pi {\Bbb N}_0$. $\blacksquare$

\section{$2$--dimensional case}\indent

In this section we will be concerned with properties of smooth cocycles
over ergodic rotations on the $2$--dimensional torus
with values in $SU(2)$. 
By $\tor^2$ we will mean the group ${\Bbb R}^2/{\Bbb Z}^2$. 
We will identify functions on $\tor^2$ with periodic of period 1 in each
coordinates functions  on ${\Bbb R}^2$.
Suppose that
$T(x_1,x_2)=(x_1+\alpha,x_2+\beta)$ is an ergodic rotation on $\tor^2$. 
Let $\varphi:\tor^2\rightarrow SU(2)$ be a $C^1$--cocycle over the
rotation $T$. Analysis similar to that in Section 2 shows that
there exists  $\psi_i\in L^2(\tor^2,\su)$, $i=1,2$ such that
\[\frac{1}{n}\frac{\partial}{\partial x_i}\varphi^{(n)}(\varphi^{(n)})^{-1}
\rightarrow\psi_i\;\;\mbox{ in }\;\;
L^2(\tor^2,\su).\]
Moreover, $\|\psi_i\|$ is a $\lambda\otimes\lambda$--a.e.\ constant function and
$\varphi(\bar{x})\psi_i(T\bar{x})\varphi(\bar{x})^{-1}=\psi_i(\bar{x})$
for $\lambda\otimes\lambda$--a.e.\ $\bar{x}\in\tor\times\tor$ for $i=1,2$.

\begin{df}
{\em The pair 
\[(\|\psi_1\|,\|\psi_2\|)=\lim_{n\rightarrow\infty}\frac{1}{n}
(\|\frac{\partial}{\partial x_1}\varphi^{(n)}(\varphi^{(n)})^{-1}\|_{L^1},
\|\frac{\partial}{\partial x_2}\varphi^{(n)}(\varphi^{(n)})^{-1}\|_{L^1})\] 
will be called the} degree {\em
of the cocycle $\varphi:\tor^2\rightarrow SU(2)$ and denoted by $d(\varphi)$.}
\end{df}

Similarly, one can prove the following

\begin{theo}
If $d(\varphi)\neq 0$, then $\varphi$ is cohomologous to a diagonal
cocycle $\tor^2\ni\bar{x}\mapsto\left[
\begin{array}{cc}
\gamma(\bar{x}) & 0\\
0 & \overline{\gamma(\bar{x})}
\end{array}\right]\in SU(2)$, where
$\gamma:\tor^2\rightarrow\tor$ is measurable. Moreover, the skew product
$T_{\gamma}:\tor^2\times\tor\rightarrow\tor^2\times\tor$ 
is ergodic and it is mixing on the  orthocomplement of the space of
functions depending only on the first two variables. $\blacksquare$
\end{theo}

Analysis similar to that in the proof of Theorem~\ref{glowne1} gives
 
\begin{theo}\label{leb2}
Let $\varphi:\tor^2\rightarrow SU(2)$ be a $C^2$--cocycle  with
$d(\varphi)\neq 0$. Suppose that the sequence
$\{\frac{1}{n}\frac{\partial}{\partial x_i}\varphi^{(n)}
(\varphi^{(n)})^{-1}\}_{n\in{\Bbb N}}$
is uniformly convergent and
$\{\frac{\partial}{\partial x_i}(\frac{1}{n}\frac{\partial}{\partial x_i}
\varphi^{(n)}(\varphi^{(n)})^{-1})\}_{n\in{\Bbb N}}$
is bounded in $L^2(\tor^2,\su)$ for $i=1,2$. Then the Lebesgue component in the
spectrum of $T_{\varphi}$ has countable multiplicity. $\blacksquare$
\end{theo}

By $BV^{{\cal R}}(\tor^2,SU(2))$ we mean the set of all measurable functions 
$f:\tor^2\rightarrow SU(2)$ 
such that 
\begin{itemize}
\item the functions $f(x,\,\cdot\,),f(\,\cdot\, ,x):\tor\rightarrow SU(2)$
are of bounded variation for any $x\in\tor$; 
\item the functions $\frac{\partial}{\partial x_1}f(f)^{-1}$,
$\frac{\partial}{\partial x_2}f(f)^{-1}:\tor^2\rightarrow \su$ are Riemann
integrable for $i=1,2$.
\end{itemize}
Then we immediately get the following

\begin{lem}\label{glas}
Let $\varphi:\tor^2\rightarrow SU(2)$ be a $C^2$--cocycle. Suppose that
$\varphi$ is cohomologous to a diagonal cocycle with a transfer
function in $BV^{{\cal R}}(\tor^2,SU(2))$.
Then the sequence
$\{\frac{1}{n}\frac{\partial}{\partial x_i}\varphi^{(n)}(\varphi^{(n)})^{-1}\}_{n\in{\Bbb N}}$
is uniformly convergent and
$\{\frac{\partial}{\partial x_i}(\frac{1}{n}\frac{\partial}{\partial x_i}\varphi^{(n)}(\varphi^{(n)})^{-1})\}_{n\in{\Bbb N}}$
is uniformly bounded for $i=1,2$. $\blacksquare$
\end{lem}

It is easy to check that if $\varphi$ is cohomologous to a diagonal
cocycle via a $C^1$ transfer function, then
$d(\varphi)\in 2\pi({\Bbb N}_0\times{\Bbb N}_0)$. However, in the next
section we show that for every ergodic rotation
$T(x_1,x_2)=(x_1+\alpha,x_2+\beta)$ there exists a smooth cocycle whose
degree is equal to $2\pi(|\beta|,|\alpha|)$.

\section{Cocycles over flows}\indent

Let $\omega$ be an irrational number. By $S:{\Bbb
R}\times\tor^2\rightarrow\tor$ we mean the ergodic flow defined by
\begin{equation}\label{flow}
S_t(x_1,x_2)=(x_1+t\omega,x_2+t). 
\end{equation}
Let  $\varphi:{\Bbb
R}\times\tor^2\rightarrow SU(2)$ be a smooth cocycle over $S$, i.e.\
\[\varphi_{t+s}(\bar{x})=\varphi_{t}(\bar{x})\varphi_{s}(S_t\bar{x})\]
for all $t,s\in{\Bbb R}$, $\bar{x}\in\tor^2$ or equivalently, $\varphi$ is
the fundamental matrix solution for a linear  differential system
\[\frac{d}{dt}y(t)=y(t)A(S_t\bar{x}),\]
where $A:\tor^2\rightarrow\su$, i.e.\ $\varphi$ satisfies
\[\left\{
\begin{array}{rcl}
\frac{d}{dt}\varphi_{t}(\bar{x}) & = & \varphi_{t}(\bar{x})A(S_t\bar{x})\\
\varphi_{0}(\bar{x}) & = & \mbox{Id}.
\end{array}
\right.\]
Then
\[\frac{\partial}{\partial x_i}\varphi_{t+s}(\bar{x})
\varphi_{t+s}(\bar{x})^{-1}=
\frac{\partial}{\partial x_i}\varphi_{t}(\bar{x})\varphi_{t}(\bar{x})^{-1}+
\ad_{\varphi_{t}(\bar{x})}\frac{\partial}{\partial x_i}\varphi_{s}(S_t\bar{x})
\varphi_{s}(S_t\bar{x})^{-1}.\]
 Hence
\[\|\frac{\partial}{\partial x_i}\varphi_{t+s}
(\varphi_{t+s})^{-1}\|_{L^1}\leq
\|\frac{\partial}{\partial x_i}\varphi_{t}(\varphi_{t})^{-1}\|_{L^1}+
\|\frac{\partial}{\partial x_i}\varphi_{s}(\varphi_{s})^{-1}\|_{L^1}\]
It follows that the limit
\[\lim_{t\rightarrow\infty}\frac{1}{|t|}
\|\frac{\partial}{\partial x_i}\varphi_{t}(\varphi_{t})^{-1}\|_{L^1}\]
exists for $i=1,2$.

\begin{df}
{\em The pair 
\[\lim_{t\rightarrow\infty}\frac{1}{|t|}
(\|\frac{\partial}{\partial x_1}\varphi_{t}(\varphi_{t})^{-1}\|_{L^1},
\|\frac{\partial}{\partial x_2}\varphi_{t}(\varphi_{t})^{-1}\|_{L^1})\]
will be called the} degree {\em
of the cocycle $\varphi:{\Bbb R}\times\tor^2\rightarrow SU(2)$ and denoted
by $d(\varphi)$.} 
\end{df}

For a given cocycle $\varphi:{\Bbb R}\times\tor^2\rightarrow SU(2)$ over
the flow $S$, by $\hp:\tor\rightarrow SU(2)$ we will mean the cocycle over the
rotation $Tx=x+\omega$ defined by $\hp(x)=\varphi_1(x,0)$. Then
$\hp^{(n)}(x)=\varphi_n(x,0)$. 

\begin{lem}
$d(\varphi)=(1,|\omega|)d(\hp)$.
\end{lem}

\pf
First observe that
\[\varphi_{x_2}(x_1-x_2\omega,0)\varphi_n(x_1,x_2)=
\varphi_{n+x_2}(x_1-x_2\omega,0)=
\varphi_n(x_1-x_2\omega,0)\varphi_{x_2}(x_1-x_2\omega+n\omega,0).\]
Hence 
\[\varphi_n(x_1,x_2)=\varphi_{x_2}(x_1-x_2\omega,0)^{-1}
\hp^{(n)}(x_1-x_2\omega)\varphi_{x_2}(x_1-x_2\omega+n\omega,0)\]
for all $x_1,x_2\in{\Bbb R}$ and $n\in{\Bbb N}$. Fix $(x_1,x_2)\in
[0,1]\times[0,1]$. Then
\begin{eqnarray*}
\lefteqn{
\frac{\partial}{\partial x_1}\varphi_n(x_1,x_2)\varphi_n(x_1,x_2)^{-1} 
=-\varphi_{x_2}(x_1-x_2\omega,0)^{-1}\frac{\partial}{\partial x_1}
\varphi_{x_2}(x_1-x_2\omega,0)}\\
& &\mbox{}+\ad_{\varphi_{x_2}(x_1-x_2\omega,0)^{-1}}(D\hp^{(n)}(x_1-x_2\omega)
\hp^{(n)}(x_1-x_2\omega)^{-1})\\
 & &\mbox{}+\ad_{\varphi_{x_2}(x_1-x_2\omega,0)^{-1}\hp^{(n)}(x_1-x_2\omega)}
(\frac{\partial}{\partial x_1}\varphi_{x_2}(x_1-x_2\omega+n\omega,0)
\varphi_{x_2}(x_1-x_2\omega+n\omega,0)^{-1}).
\end{eqnarray*}
It follows that
\begin{eqnarray*}
\lefteqn{|\|\frac{\partial}{\partial x_1}\varphi_n(\varphi_n)^{-1} \|_{L^1}-
\|D\hp^{(n)}(\hp^{(n)})^{-1}\|_{L^1}|} & & \\ & = &
 |\|\frac{\partial}{\partial x_1}\varphi_n(\varphi_n)^{-1} \|_{L^1}-
\int_0^1\int_0^1\|D\hp^{(n)}(x_1-x_2\omega)\hp^{(n)}(x_1-x_2\omega)^{-1}
\|dx_1dx_2|\\
& \leq & 2\int_0^1\int_0^1\|\frac{\partial}{\partial x_1}
\varphi_{x_2}(x_1-x_2\omega,0)
\varphi_{x_2}(x_1-x_2\omega,0)^{-1}
\|dx_1dx_2.
\end{eqnarray*}
Therefore
\[\lim_{n\rightarrow\infty}\frac{1}{n}
\|\frac{\partial}{\partial x_1}\varphi_n(\varphi_n)^{-1} \|_{L^1}=
\lim_{n\rightarrow\infty}\frac{1}{n}
\|D\hp^{(n)}(\hp^{(n)})^{-1}\|_{L^1}=d(\hp).\]
Similarly,
\begin{eqnarray*}
\lefteqn{
\frac{\partial}{\partial x_2}\varphi_n(x_1,x_2)\varphi_n(x_1,x_2)^{-1} 
=-\varphi_{x_2}(x_1-x_2\omega,0)^{-1}\frac{\partial}{\partial t}
\varphi_{x_2}(x_1-x_2\omega,0)}\\
 & & \mbox{}+\omega\,\varphi_{x_2}(x_1-x_2\omega,0)^{-1}\frac{\partial}{\partial x_1}
\varphi_{x_2}(x_1-x_2\omega,0)\\ 
 & & \mbox{}-\omega\ad_{\varphi_{x_2}(x_1-x_2\omega,0)^{-1}}(D\hp^{(n)}(x_1-x_2\omega)
\hp^{(n)}(x_1-x_2\omega)^{-1})\\
 & & \mbox{}+\ad_{\varphi_{x_2}(x_1-x_2\omega,0)^{-1}\hp^{(n)}(x_1-x_2\omega)}
(\frac{\partial}{\partial t}\varphi_{x_2}(x_1-x_2\omega+n\omega,0)
\varphi_{x_2}(x_1-x_2\omega+n\omega,0)^{-1})\\
 & & \mbox{}-\omega\ad_{\varphi_{x_2}(x_1-x_2\omega,0)^{-1}\hp^{(n)}(x_1-x_2\omega)}
(\frac{\partial}{\partial x_1}\varphi_{x_2}(x_1-x_2\omega+n\omega,0)
\varphi_{x_2}(x_1-x_2\omega+n\omega,0)^{-1}).
\end{eqnarray*}
It follows that
\begin{eqnarray*}
\lefteqn{|\|\frac{\partial}{\partial x_2}\varphi_n(\varphi_n)^{-1} \|_{L^1}-
|\omega|\|D\hp^{(n)}(\hp^{(n)})^{-1}\|_{L^1}|} & & \\ 
& \leq & 2\int_0^1\int_0^1\|\frac{\partial}{\partial t}
\varphi_{x_2}(x_1-x_2\omega,0)
\varphi_{x_2}(x_1-x_2\omega,0)^{-1}
\|dx_1dx_2\\ 
& &\mbox{}+2|\omega|\int_0^1\int_0^1\|\frac{\partial}{\partial x_1}
\varphi_{x_2}(x_1-x_2\omega,0)
\varphi_{x_2}(x_1-x_2\omega,0)^{-1}
\|dx_1dx_2.
\end{eqnarray*}
Therefore
\[\lim_{n\rightarrow\infty}\frac{1}{n}
\|\frac{\partial}{\partial x_2}\varphi_n(\varphi_n)^{-1} \|_{L^1}=
|\omega|\lim_{n\rightarrow\infty}\frac{1}{n}
\|D\hp^{(n)}(\hp^{(n)})^{-1}\|_{L^1}=|\omega| d(\hp),\]
and the proof is complete. $\blacksquare$

\begin{lem}
For any $C^2$--cocycle $\psi:\tor\rightarrow SU(2)$ over the rotation $T$
there exists a $C^2$--cocycle $\varphi:{\Bbb R}\times\tor^2\rightarrow
SU(2)$ over the flow $S$ such that $\hp=\psi$.  
\end{lem}

\pf
Since the fundamental group of $SU(2)$ is trivial, we can choose a
$C^2$--homotopy $\psi:[0,1]\times\tor\rightarrow SU(2)$ such that
\[\psi(t,x)=\left\{
\begin{array}{rcl}
\mbox{Id} & \mbox{for} & t\in[0,1/4]\\
\psi(x) & \mbox{for} & t\in[3/4,1].
\end{array}
\right.\]
By $\psi:{\Bbb R}\times\tor\rightarrow SU(2)$ we mean the $C^2$--function
determined by
\[\psi(n+t,x)=\psi^{(n)}(x)\psi(t,x+n\omega)\]
for any $t\in[0,1]$ and $n\in{\Bbb Z}$.
Then it is easy to check that 
\begin{equation}\label{nat}
\psi(n+t,x)=\psi^{(n)}(x)\psi(t,x+n\omega)
\end{equation}
for any $t\in{\Bbb R}$ and $n\in{\Bbb Z}$.
Let $\varphi:{\Bbb R}\times{\Bbb R}^2\rightarrow SU(2)$ be defined by
\[\varphi_t(x_1,x_2)=\psi(x_2,x_1-x_2\omega)^{-1}\psi(t+x_2,x_1-x_2\omega).\]
It is easy to see that $\varphi_t(x_1+1,x_2)=\varphi_t(x_1,x_2)$ and
$\varphi_t(x_1,x_2+1)=\varphi_t(x_1,x_2)$, by (\ref{nat}).
Then $\varphi:{\Bbb R}\times\tor^2\rightarrow SU(2)$ is a $C^2$--function
and
\begin{eqnarray*}
\varphi_{t+s}(\bar{x}) & = & \psi(x_2,x_1-x_2\omega)^{-1}\psi(t+s+x_2,x_1-x_2\omega)\\
 & = & \psi(x_2,x_1-x_2\omega)^{-1}\psi(t+x_2,x_1-x_2\omega)\\
 & &\psi(x_2+t,(x_1+t\omega)-(x_2+t)\omega)^{-1}\psi(s+(x_2+t),(x_1+t\omega)-(x_2+t)\omega)\\
 & = & \varphi_{t}(\bar{x})\varphi_{s}(S_t\bar{x}).
\end{eqnarray*}
Moreover,
\[\hp(x)=\varphi_1(x,0)=\psi(0,x)^{-1}\psi(1,x)=\psi(x),\]
which completes the proof. $\blacksquare$

\vspace{2ex}

Suppose that $\alpha,\beta,1$ are independent over ${\Bbb Q}$. Set
$\omega=\alpha/\beta$. 

\begin{theo}
For every ergodic rotation $T(x_1,x_2)=(x_1+\alpha,x_2+\beta)$ and for
every natural $k$  there exists a $C^2$--cocycle over $T$
whose degree is equal to $2\pi k(|\beta|,|\alpha|)$.
\end{theo}

\pf
Let $S$ denote the ergodic flow given by (\ref{flow}).
Suppose that $\varphi:{\Bbb R}\times\tor^2\rightarrow SU(2)$
is a $C^2$--cocycle over $S$ such that $d(\hp)=2\pi k$. Consider the cocycle
$\varphi_{\beta}:\tor\rightarrow SU(2)$ over the rotation
$T=S_{\beta}$. Then $\varphi_{\beta}^{(n)}=\varphi_{\beta n}$ and
\[\lim_{n\rightarrow\infty}\frac{1}{n}\|\frac{\partial}{\partial x_i}
\varphi_{\beta}^{(n)}(\varphi_{\beta}^{(n)})^{-1}\|_{L_1}=
|\beta|
\lim_{n\rightarrow\infty}\frac{1}{|\beta| n}\|\frac{\partial}{\partial x_i}
\varphi_{\beta n}(\varphi_{\beta n})^{-1}\|_{L_1}.\]
It follows that
\[d(\varphi_{\beta})=|\beta| d(\varphi)=|\beta|(1,|\omega|)\,d(\hp)
=(|\beta|,|\alpha|)\,d(\hp),\]
which proves the theorem. $\blacksquare$

\vspace{2ex}

Suppose that $\beta\in(0,1)$.
Let $\varphi:{\Bbb R}\times\tor^2\rightarrow SU(2)$
be a $C^2$--cocycle over $S$ such that $\hp$ is a diagonal $C^2$--cocycle
with nonzero degree. Set $T=S_{\beta}$ and $\psi=\varphi_{\beta}$. Let
$p:\tor^2\rightarrow SU(2)$ be a $BV^{{\cal R}}$--function such that
\[p(x_1,x_2)=\varphi_{x_2}(x_1-x_2\omega,0)^{-1}\]
for $(x_1,x_2)\in{\Bbb R}\times[0,1)$. Then
\[p(T(x_1,x_2))=\left\{
\begin{array}{rcl}
\varphi_{x_2+\beta}(x_1-x_2\omega,0)^{-1} & \mbox{for} & x_2\in[0,1-\beta)\\
\varphi_{x_2+\beta-1}(x_1-(x_2-1)\omega,0)^{-1} & \mbox{for} & x_2\in[1-\beta,1).
\end{array}
\right.\]
Moreover,
\[\varphi_{x_2+\beta}(x_1-x_2\omega,0)=\varphi_{x_2}(x_1-x_2\omega,0)
\varphi_{\beta}(x_1,x_2)\]
and
\begin{eqnarray*}
\varphi_{x_2+\beta-1}(x_1-(x_2-1)\omega,0) & = & 
\varphi_{-1}(x_1-(x_2-1)\omega,0)
\varphi_{x_2+\beta}(x_1-x_2\omega,0)\\
 & = & 
\varphi_1(x_1-x_2\omega,0)^{-1}
\varphi_{x_2+\beta}(x_1-x_2\omega,0).
\end{eqnarray*}
It follows that $p(\bar{x})\delta(\bar{x})p(T\bar{x})^{-1}=\psi(\bar{x})$,
where $\delta:\tor^2\rightarrow SU(2)$ is the diagonal $BV^{{\cal
R}}$--cocycle given by
\[\delta(x_1,x_2)=\left\{
\begin{array}{rcl}
\mbox{Id} & \mbox{for} & x_2\in[0,1-\beta)\\
\hp(x_1-x_2\omega) & \mbox{for} & x_2\in[1-\beta,1).
\end{array}
\right.\]

\begin{lem}\label{cob}
Let $\phi:\tor^2\rightarrow\tor$ be a cocycle over the rotation 
$T(x_1,x_2)=(x_1+\alpha,x_2+\beta)$. Suppose that  
$\phi|\tor\times[0,\gamma)$, $\phi|\tor\times[\gamma,1)$ are
$C^1$--functions, where $\gamma$ is irrational. If 
$d(\phi(\,\cdot\, ,0))\neq d(\phi(\,\cdot\, ,\gamma))$, then $\phi$
is not a coboundary.  
\end{lem}

\pf Set $I_1=[0,\gamma)$, $I_2=[\gamma,1)$, $a_1=d(\phi(\,\cdot\, ,0))$ and
$a_2=d(\phi(\,\cdot\, ,\gamma))$. Then there exists a function
$\ts:\tor^2\rightarrow{\Bbb R}$  such that $\ts|\tor\times I_j$
is of class $C^1$ for $j=1,2$ and
$\phi(x_1,x_2)=\exp 2\pi i(\ts(x_1,x_2)+a_jx_1)$ for any 
$(x_1,x_2)\in \tor\times I_j$.

Clearly, it suffices to show that
\[\int_{\tor^2} \phi^{(n)}(x_1,x_2)dx_1dx_2\rightarrow 0.\]
Next note that
\[\phi^{(n)}(x_1,x_2)=\exp 2\pi i(\ts^{(n)}(x_1,x_2)+
(a_1S_1^n(x_2)+a_2S_2^n(x_2))x_1+c_n(x_2)),\]
where $S_i^n(x)=\sum_{k=0}^{n-1}\jed_{I_i}(x+k\beta)$ and
$c_n(x)=\sum_{k=0}^{n-1}k\alpha(a_1\jed_{I_1}+a_2\jed_{I_2})(x+k\beta)$.
 Since the rotation by $\beta$ is uniquely ergodic, 
\[\frac{1}{n}(a_1S_1^n+a_2S_2^n)\rightarrow a_1\gamma+a_2(1-\gamma)\]
uniformly.
Since $a_1\neq a_2$ and $\gamma$ is irrational, there exists $S>0$ and
$n_0\in{\Bbb N}$ such that $|a_1S_1^n(x)+a_2S_2^n(x)|\geq nS$ for all
$x\in\tor$ and $n\geq n_0$. Applying integration by parts, we get  
\begin{eqnarray*}
\lefteqn{|\int_{\tor^2} \phi^{(n)}(x_1,x_2)dx_1dx_2|}\\ & \leq &
\int_0^1|\int_0^1e^{2\pi i\ts^{(n)}(x_1,x_2)+(a_1S_1^n(x_2)+a_2S_2^n(x_2))
x_1}dx_1|dx_2\\
 & = &
\int_0^1
\frac{1}{2\pi |a_1S_1^n(x_2)+a_2S_2^n(x_2)|}
|\int_0^1e^{2\pi i\ts^{(n)}(x_1,x_2)}
de^{2\pi i(a_1S_1^n(x_2)+a_2S_2^n(x_2))x_1}|dx_2\\
 & = &
\int_0^1
\frac{1}{2\pi |a_1S_1^n(x_2)+a_2S_2^n(x_2)|}
|\int_0^1e^{2\pi i(a_1S_1^n(x_2)+a_2S_2^n(x_2))x_1}
de^{2\pi i\ts^{(n)}(x_1,x_2)}|dx_2\\
 & \leq &
\int_0^1
\frac{1}{nS}
|\int_0^1e^{2\pi i\ts^{(n)}(x_1,x_2)+
(a_1S_1^n(x_2)+a_2S_2^n(x_2))x_1}
\frac{\partial}{\partial x_1}\ts^{(n)}(x_1,x_2)dx_1|dx_2\\
 & \leq &
\frac{1}{nS}
\int_{\tor^2}|\frac{\partial}{\partial x_1}\ts^{(n)}(x_1,x_2)|dx_1dx_2.
\end{eqnarray*}
Since $\frac{\partial}{\partial x_1}\ts\in L^1(\tor^2,{\Bbb C})$,
\[\frac{1}{n}\frac{\partial}{\partial x_1}\ts^{(n)}\rightarrow
\int_{\tor^2}\frac{\partial}{\partial x_1}\ts(x_1,x_2)dx_1dx_2=0\]
in $L^1(\tor^2,{\Bbb C})$, by the Birkhoff ergodic theorem, and the proof in complete.
$\blacksquare$ 

\vspace{2ex}

This leads to the following conclusion.

\begin{cor} 
For every ergodic rotation $T$ on $\tor^2$ there exists a $C^2$--cocycle
$\psi$ with nonzero degree such that the Lebesgue component in the spectrum of
$T_{\psi}$ has countable multiplicity and $\psi$ is not cohomologous 
to any diagonal $C^1$--cocycle. 
\end{cor}

\pf
Let $\check{\varphi}:\tor\rightarrow\tor$ be a $C^2$--function with
nonzero topological degree. Let 
$\varphi:{\Bbb R}\times\tor^2\rightarrow SU(2)$
be a $C^2$--cocycle over $S$ such that 
$\hp=\left[\begin{array}{cc}
\check{\varphi} & 0\\
0 & (\check{\varphi})^{-1}
\end{array}\right]$.
Define $\psi=\varphi_{\beta}$.
Then $d(\psi)=2\pi(|\beta|,|\alpha|)|d(\check{\varphi})|\neq 0$.
Moreover, $\psi$ and the
diagonal cocycle  $\delta:\tor^2\rightarrow SU(2)$ given by
\[\delta(x_1,x_2)=\left\{
\begin{array}{rcl}
\mbox{Id} & \mbox{for} & x_2\in[0,1-\beta)\\
\hp(x_1-x_2\omega) & \mbox{for} & x_2\in[1-\beta,1)
\end{array}
\right.\]
are cohomologous with a transfer function in $BV^{{\cal R}}(\tor^2,SU(2))$.
Applying Theorem~\ref{leb2} and Lemma~\ref{glas}, we get the first part of
our claim. 

Next suppose that $\psi$ is cohomologous to a diagonal $C^1$--cocycle.
Then it is easy to see that the cocycle $\eta:\tor^2\rightarrow\tor$ 
given by 
\[\eta(x_1,x_2)=\left\{
\begin{array}{rcl}
\mbox{Id} & \mbox{for} & x_2\in[0,1-\beta)\\
\check{\varphi}(x_1-x_2\omega) & \mbox{for} & x_2\in[1-\beta,1).
\end{array}
\right.\]
is cohomologous to a $C^1$--cocycle $g:\tor^2\rightarrow\tor$. 
Applying Lemma~\ref{cob} for $\phi=\eta\,g^{-1}$ and $\gamma=1-\beta$ we
find that $\eta\,g^{-1}$ is 
not a coboundary, which completes the proof. $\blacksquare$

\noindent 
Faculty of Mathematics and Computer Science,\\
Nicholas Copernicus University\\ 
ul. Chopina 12/18\\ 
87-100 Toru\'n, Poland \\
E-mail: fraczek@mat.uni.torun.pl\\
\hspace*{8ex}Krzysztof.Fraczek@esi.ac.at

\end{document}